\documentclass[reqno]{amsart}

\usepackage{amsmath,amssymb,amsthm}
\usepackage{graphicx,tikz}
\newtheorem{theorem}{Theorem}
\newtheorem*{thm}{Theorem}
\newtheorem{proposition}{Proposition}
\newtheorem*{prop}{Proposition}
\newtheorem*{corollary}{Corollary}

\theoremstyle{definition}

\theoremstyle{remark}

\DeclareMathOperator{\diam}{diam}

\begin{document}

\title[]{Curvature on Graphs via equilibrium measures}
\subjclass[2010]{05C99, 31C20, 91A80.} 
\keywords{Graph, Curvature, Bonnet-Myers, Lichnerowicz, Minimax Theorem.}
\thanks{S.S. is supported by the NSF (DMS-2123224) and the Alfred P. Sloan Foundation.}

\author[]{Stefan Steinerberger}
\address{Department of Mathematics, University of Washington, Seattle, WA 98195, USA}
\email{steinerb@uw.edu}

\begin{abstract} We introduce a notion of curvature on finite, combinatorial graphs. It can be easily computed by solving a linear system of equations. We show that graphs with curvature bounded below by $K>0$ have diameter bounded by $\diam(G) \leq 2/K$ (a Bonnet-Myers theorem), that $\diam(G) = 2/K$ implies that $G$ has constant curvature (a Cheng theorem) and that there is a spectral gap $\lambda_1 \geq K/(2n)$ (a Lichnerowicz theorem). It is computed for several families of graphs and often coincides with Ollivier curvature or Lin-Lu-Yau curvature. The von Neumann minimax theorem features prominently in the proofs.
\end{abstract}

\maketitle

\vspace{-8pt}

\section{Introduction} 
\subsection{Introduction.} The notion of curvature is one of the cornerstones of differential geometry and geometric analysis. Starting with the work of Bakry-\'Emery \cite{bakry} there has been substantial interest in defining curvature in more abstract spaces and on graphs. While there are purely combinatorial definitions  
 \cite{higuchi, stone, woess}, many of these notions are inspired by the behavior of the Laplacian (Bakry-\'Emery curvature \cite{bakry} or Forman curvature \cite{forman}) or the behavior of optimal
transport (Lott-Villani \cite{lott}, Sturm \cite{sturm}), for example Ollivier-Ricci curvature \cite{olli0, ollivier, olli2} and the Lin-Lu-Yau curvature \cite{lly} (both are defined on edges instead of vertices). This is an active field of research, we do not aim to give a complete overview here and instead refer to
 \cite{ambrosio, bauer0, bauer, bauer2, rigidity, horn, jost, linyau, maas, villani} and references therein.

 \subsection{Definition.} We define a potential-theoretic notion of curvature by looking for a signed measure $\mu$ defined on the vertices of the graph such that
 $$ \forall v \in V \qquad \quad \sum_{w \in V} d(v,w) \mu(w) = |V|.$$
 Given $\mu$, we will interpret $\mu(v)$ as the curvature of the graph in the vertex $v \in V$.   
 \begin{center}
\begin{figure}[h!]
\begin{tikzpicture}[scale=1]
\node at (0,-3) {\includegraphics[width=0.24\textwidth]{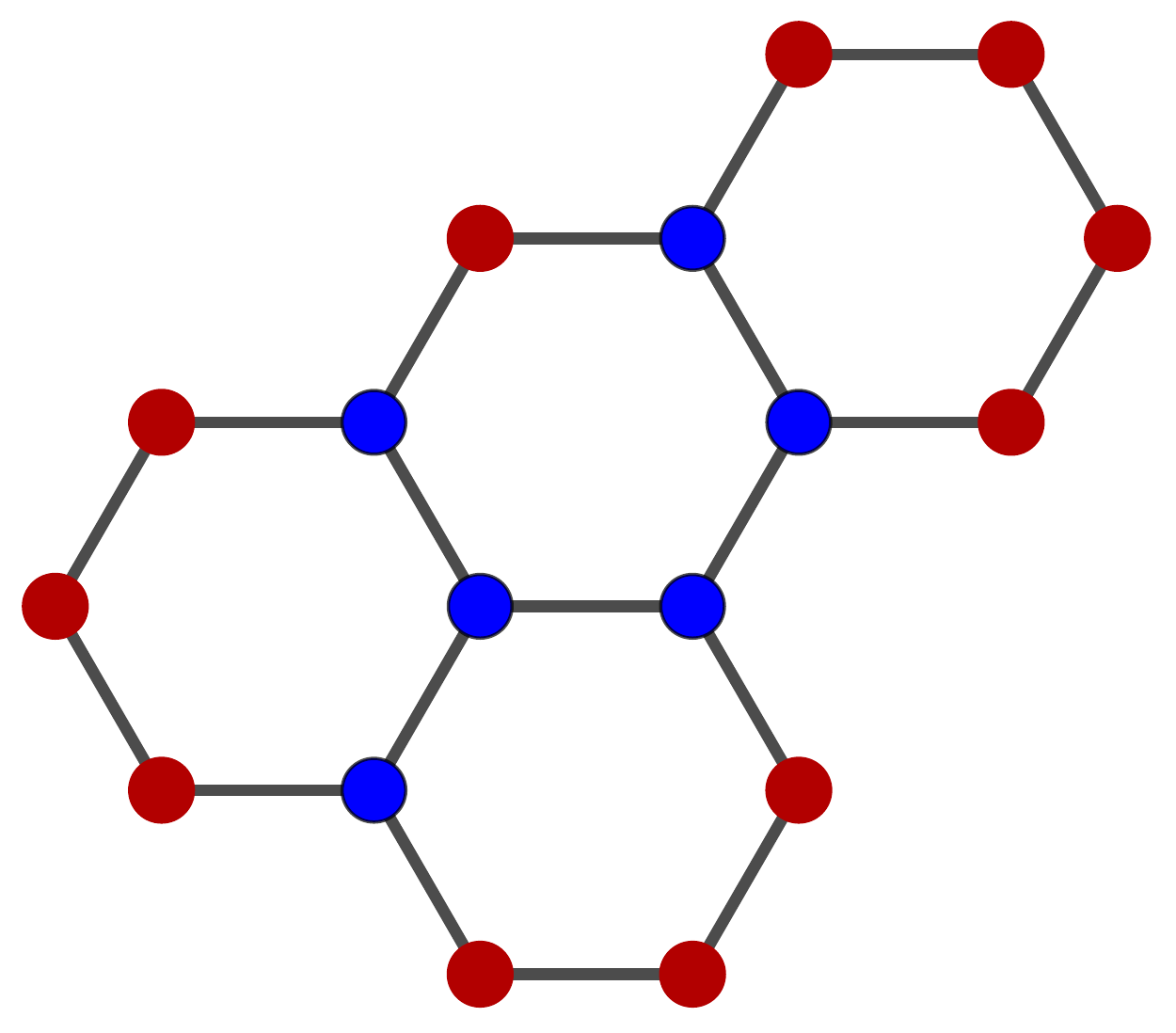}};
\node at (3.5,-3) {\includegraphics[width=0.28\textwidth]{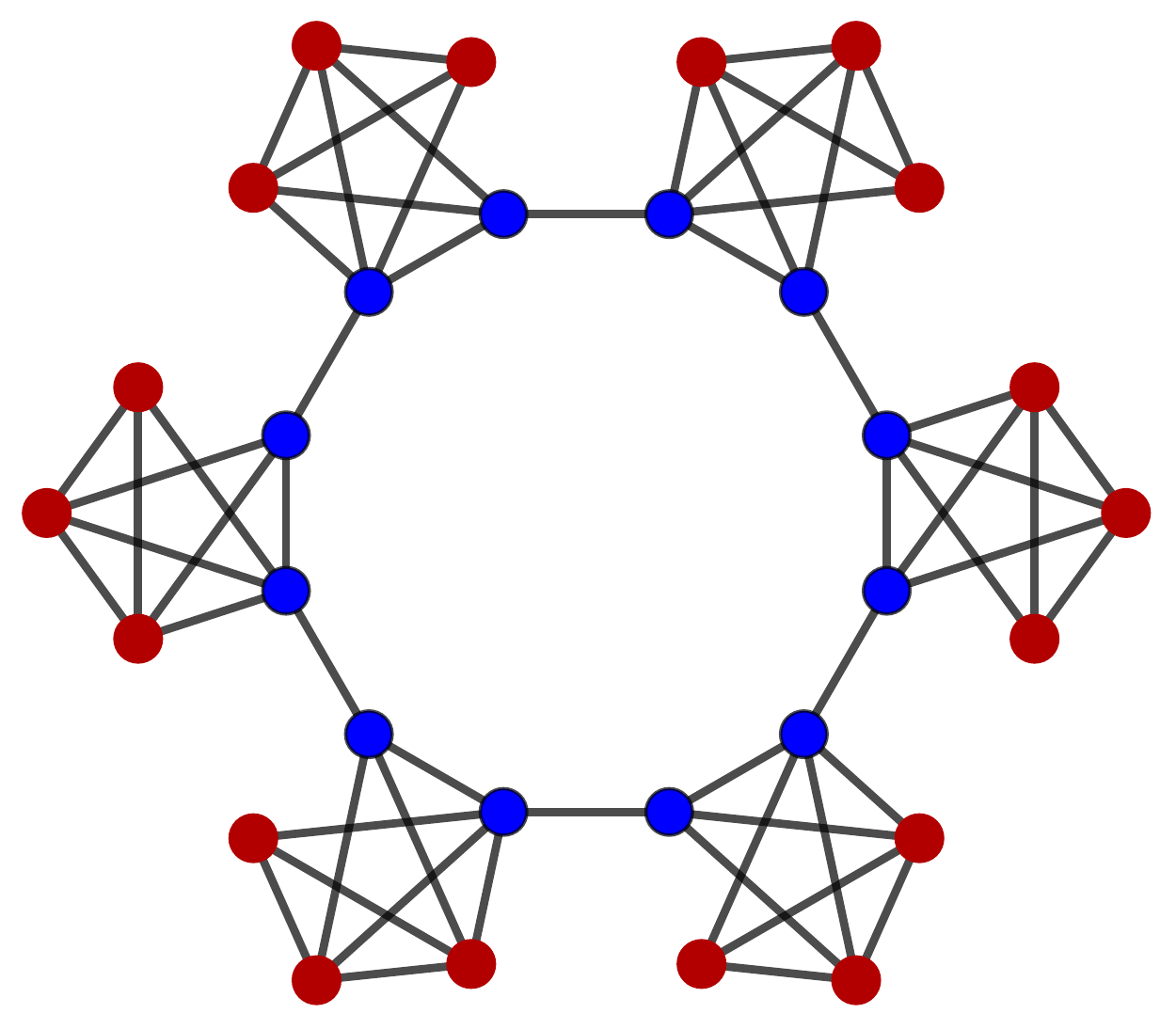}};
\node at (7,-3) {\includegraphics[width=0.22\textwidth]{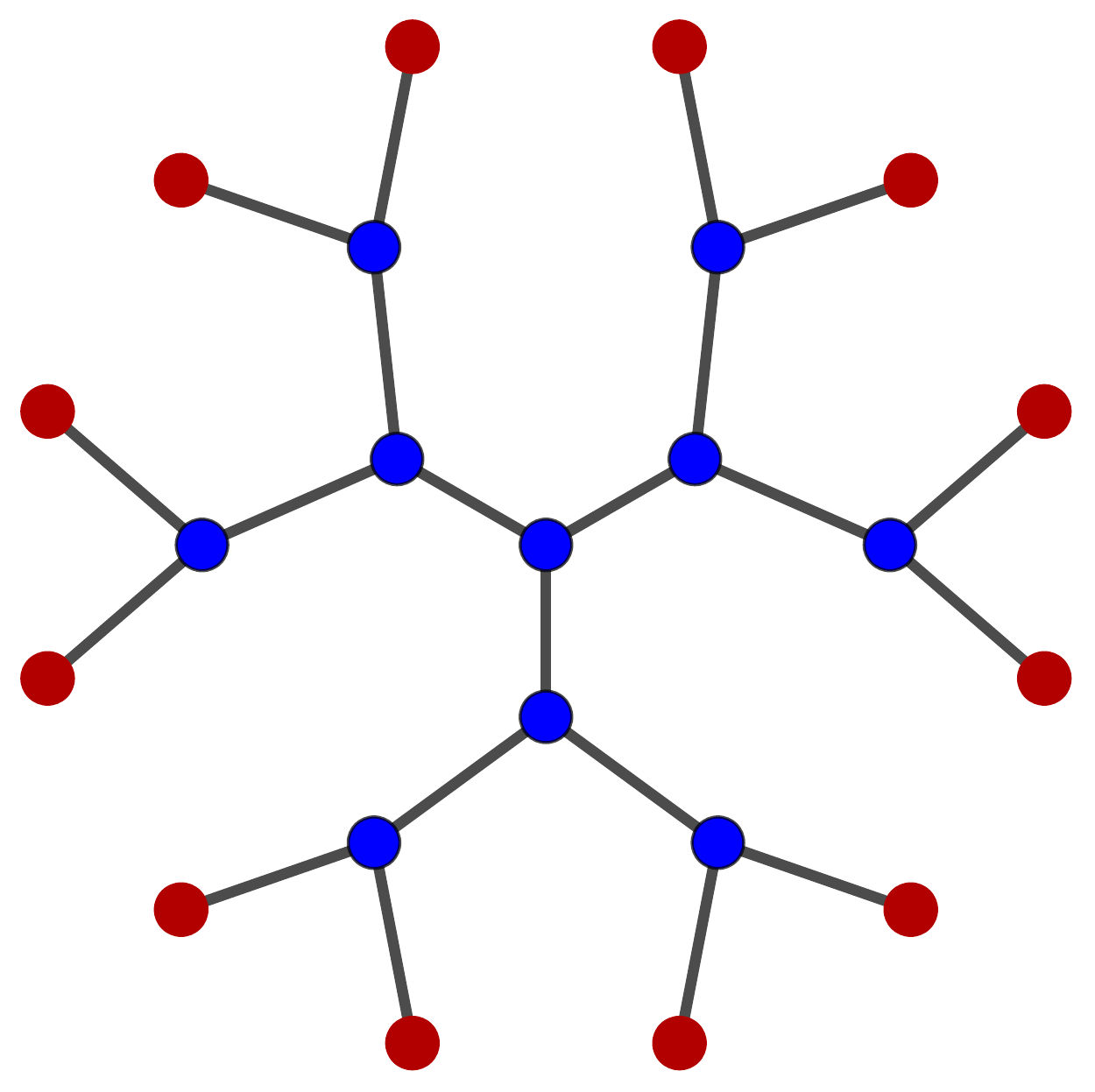}};
\end{tikzpicture}
\vspace{-10pt}
\caption{Vertices colored by curvature (red if positive, blue if negative).}
\end{figure}
\end{center}

An example of what positively and negatively curved regions look like is given in Fig. 1. We start by noting that this curvature is very easy to compute: if $G$ is a graph on $n = |V|$ vertices, then any such equilibrium measure $\mu$ corresponds to a vector of weights $w = (\mu(v_1), \mu(v_2), \dots, \mu(v_n)) \in \mathbb{R}^n$ solving the linear system
$$ D w = n \cdot \mathbf{1},~\mbox{where} \quad D =\left( d(v_i, v_j) \right)_{i,j = 1}^{n}$$
is the distance matrix and $\mathbf{1} \in \mathbb{R}^n$ is the vector containing all 1's. Finding the curvature vector $w$ therefore only requires solving a linear system of equations. There are two natural questions: does the linear system have a solution and is it unique? As it turns out, uniqueness is not much of an issue: if a graph has nonnegative curvature, meaning there exists $w \in \mathbb{R}_{\geq 0}^{n}$ such that $ Dw = n \cdot \mathbf{1}$, then the total curvature $\|w\|_{\ell^1}$ is an invariant. If there are multiple solutions, a canonical choice is one maximizing the lower curvature bound $K= \min_i w_i$. Albeit seemingly very rare (see \S 2.5), there are graphs for which $ D w = n \cdot \mathbf{1}$ has no solution: for these we can consider the Moore-Penrose pseudoinverse $w = D^{\dagger}\left(n \cdot \mathbf{1} \right)$ (which always exists). This case is dealt with in \S 2.5 where a generalized Bonnet-Myers theorem and a generalized Lichnerowicz theorem for this setting is established.
The definition was discovered coincidentally while considering unrelated problems in potential theory.
$$ \sum_{w \in V} d(v,w) \mu(w) = |V|$$
can be interpreted as an equilibrium measure: $\mu$ is a signed measure which is acting proportional to the distance and we ask for the measure to be such that the left-hand side (which can be interpreted as an effective force) is in equilibrium. The purpose of our paper is to point out that this definition has a large number of desirable and interesting properties.
\begin{enumerate}
\item There are no additional parameters that one needs to tune.
\item It satisfies a Bonnet-Myers theorem (Theorem 1), a Cheng Theorem (also Theorem 1) and a Lichnerowicz Theorem (Theorem 3). These three Theorems are sometimes considered minimal requirements for a notion of curvature to be reasonable. Moreover, it satisfies a \textit{reverse} Bonnet-Myers Inequality (Theorem 2) and, strongest of all, a Minimax Theorem (Theorem 4) which is stronger than all previous results. Indeed, we will derive all previous results as a consequence of Theorem 4.
\item This notion of curvature tends to lead to very similar (and often exactly the same) answers as the Ollivier \cite{ollivier} curvature or the Lin-Lu-Yau \cite{lly} curvature (see \S 1.3). It is not presently clear why that is the case. 
\item It is much easier to compute than curvatures based on Optimal Transport, it only requires the solution of a linear system of equations. It is also easier to compute in closed form for explicit families of graphs (see \S 1.3).
\item The linear system of equations $D \mu = |V|$ seems to have a solution in virtually all cases (in the sense that exceptions, for which a substitute theory (Theorem 5) is presented, are exceedingly rare). This is a perhaps unexpected phenomenon and interesting in its own right, see \S 2.5 and \cite{stein}.
\item Finally, in contrast to other curvatures, the cycle graph $C_n$ has constant curvature $K = n/ \left\lfloor n^2/4 \right\rfloor$, inversely proportional to diameter, which is the natural scaling one would expect from continuous considerations. \end{enumerate}

\subsection{Examples.}
We start by discussing some examples (see Fig. 2). What is somewhat remarkable is that even though our definition is quite different from Ollivier curvature and Lin-Lu-Yau curvature (both of which are defined on edges instead of vertices), we get similar or even identical results in many cases. \\

\begin{center}
\begin{figure}[h!]
\begin{tikzpicture}
\node at (-0.5,0) {\includegraphics[width=0.24\textwidth]{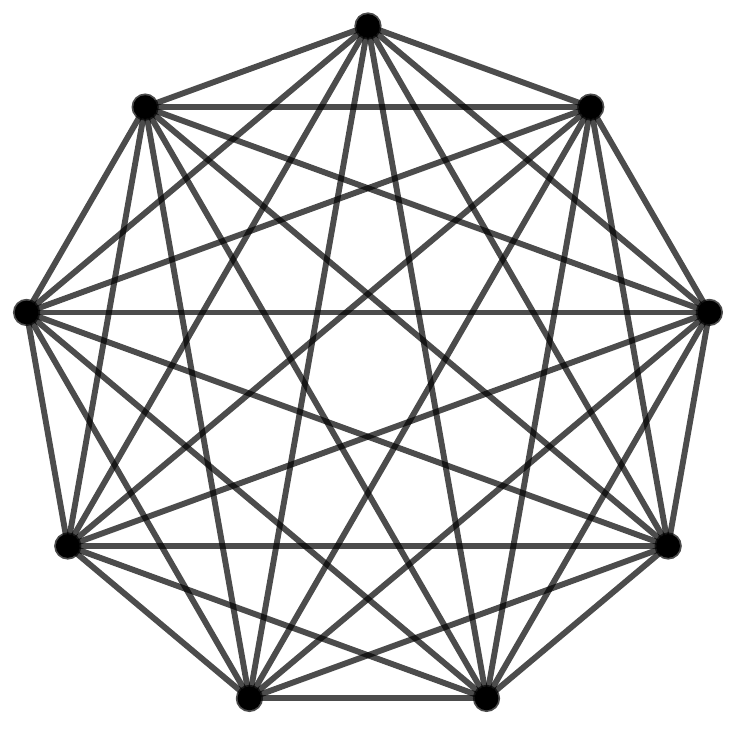}};
\node at (-0.5, -1.8) {$K(K_n) = \frac{n}{n-1}$};
\node at (3,0) {\includegraphics[width=0.2\textwidth]{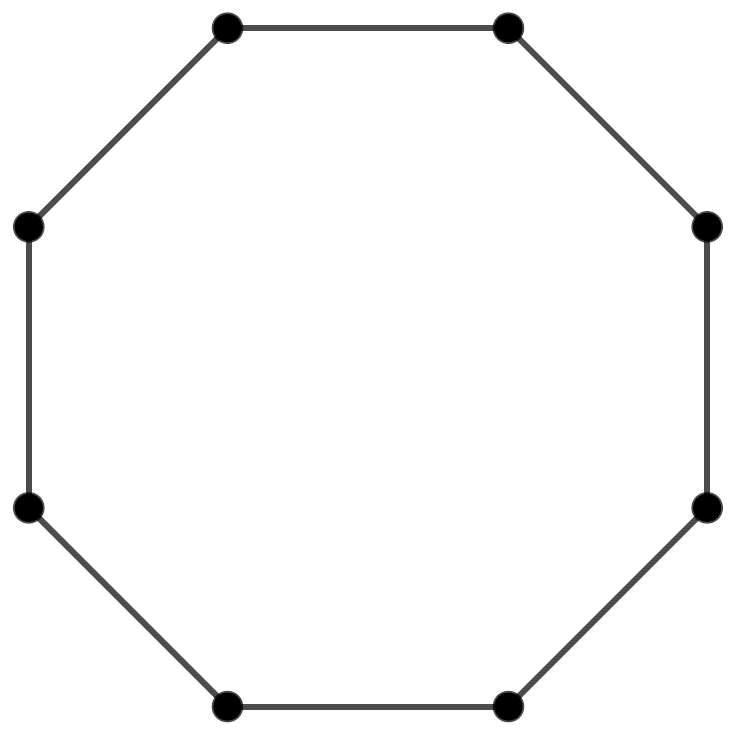}};
\node at (3, -1.8) {$K(C_n) = \frac{n}{\left\lfloor n^2/4 \right\rfloor }$};
\node at (7,0) {\includegraphics[width=0.24\textwidth]{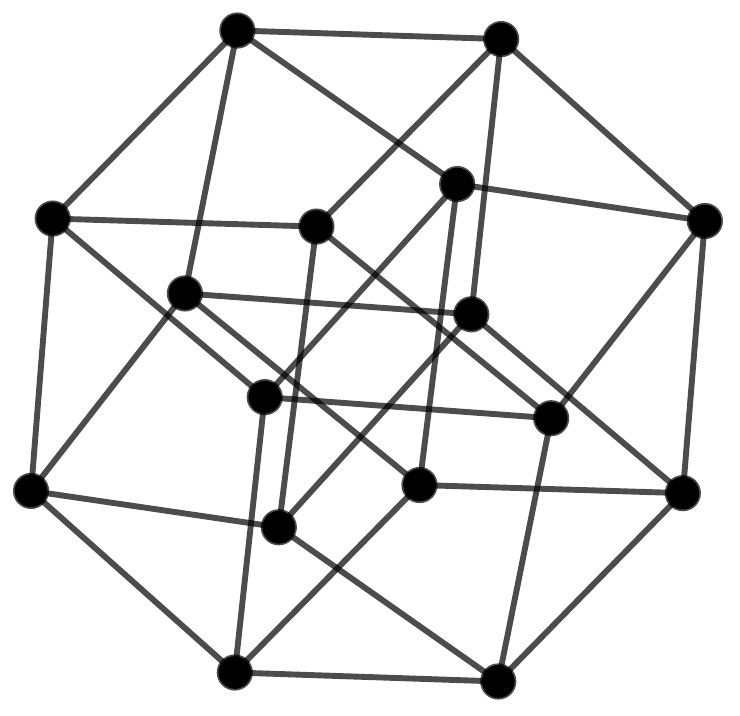}};
\node at (7, -1.8) {$K(Q_n) = \frac{2}{n} $};
\node at (-0.5,-4) {\includegraphics[width=0.24\textwidth]{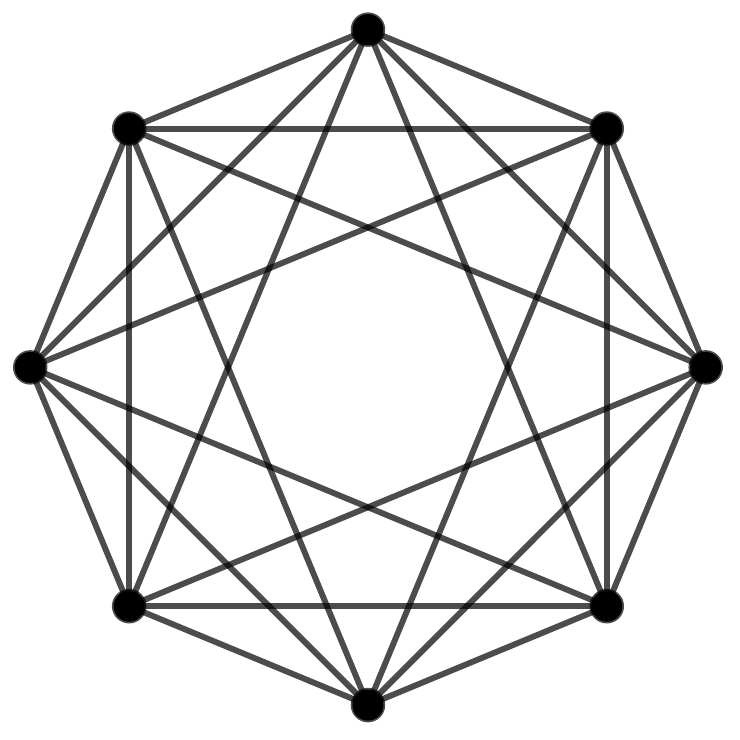}};
\node at (-0.5, -5.8) {$K(CP_n) = 1$};
\node at (3,-4) {\includegraphics[width=0.24\textwidth]{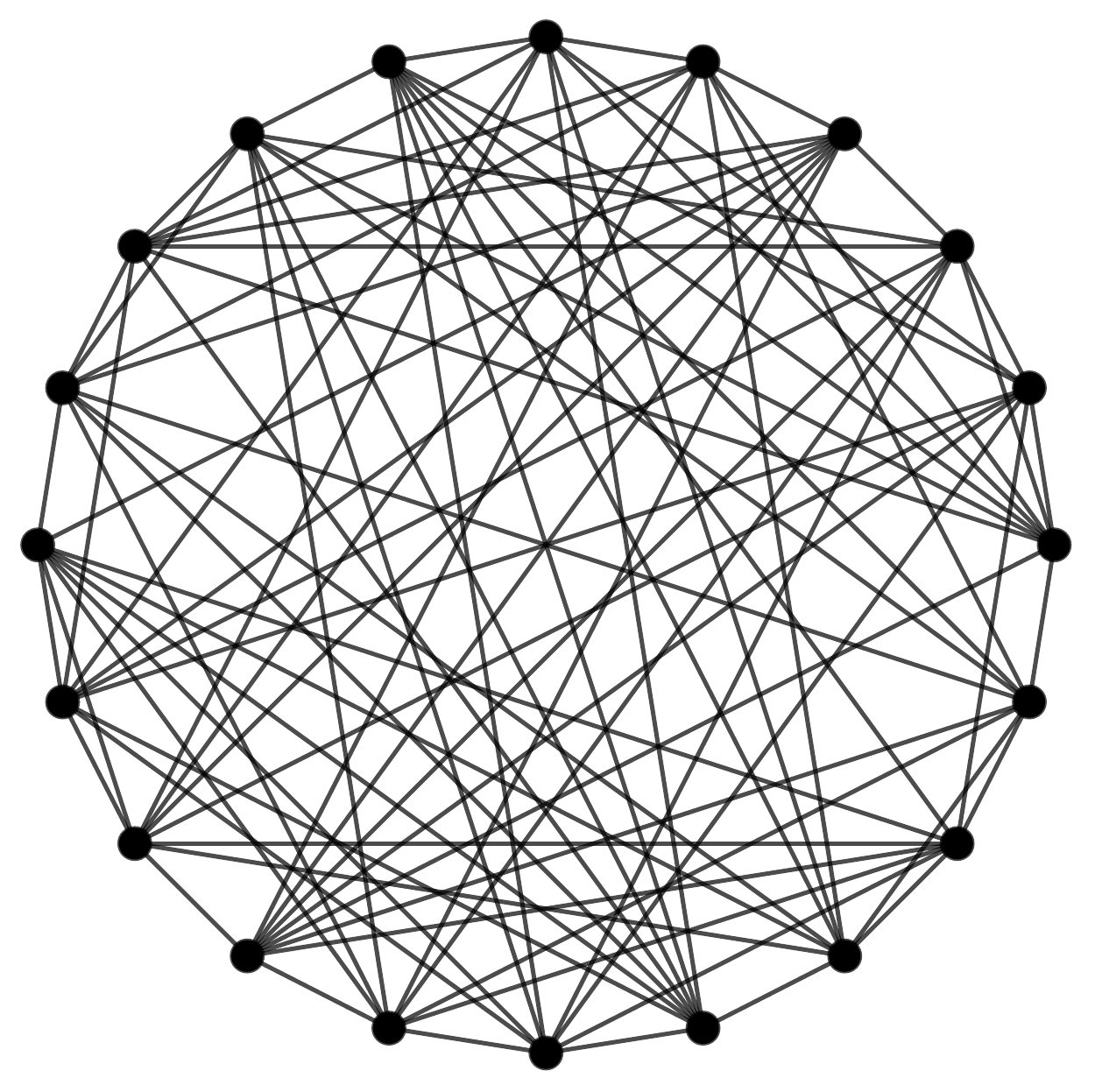}};
\node at (3, -5.8) {$K(J_{n,k}) = \frac{n}{k(n-k)}$};
\node at (7,-4) {\includegraphics[width=0.24\textwidth]{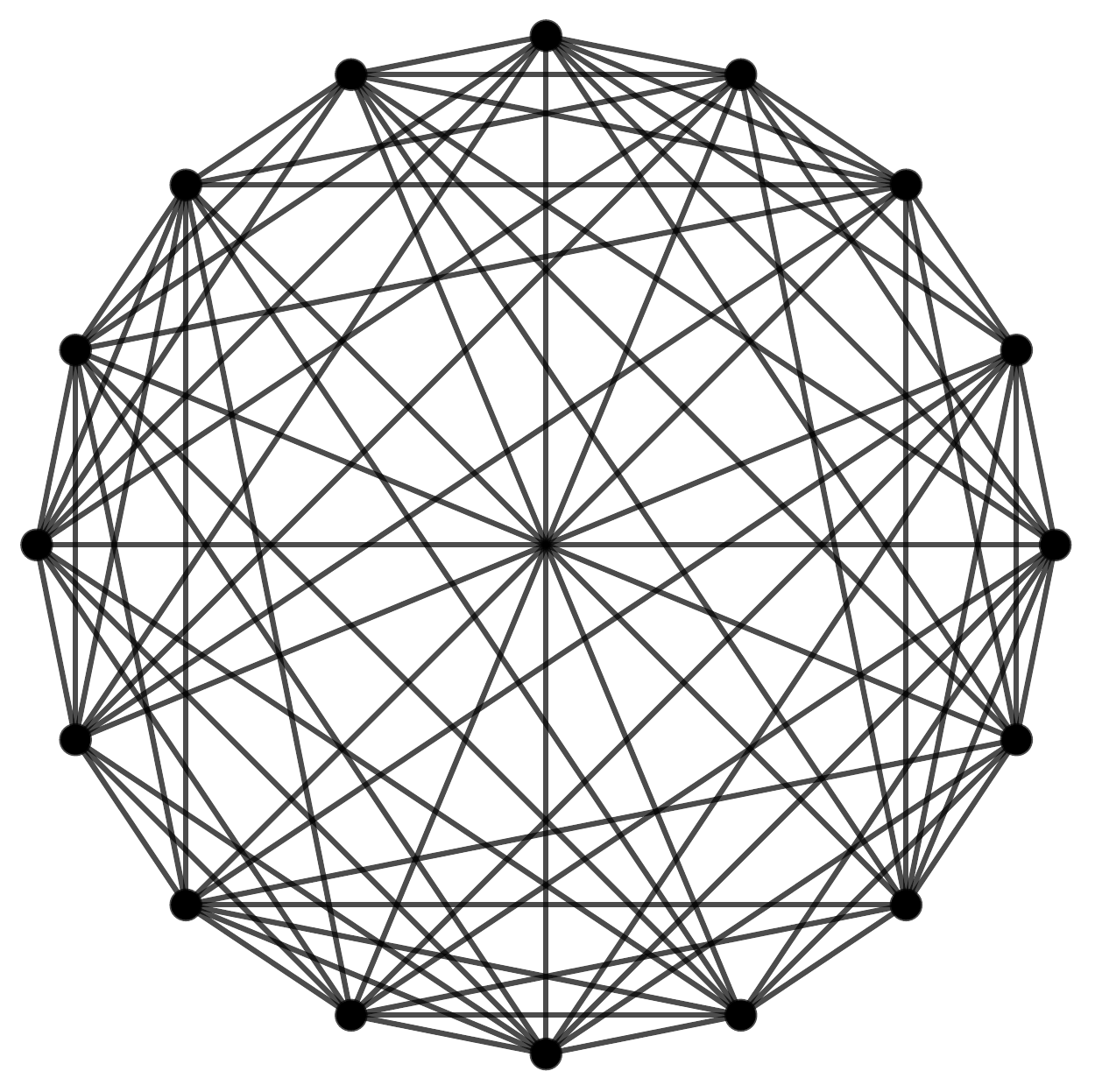}};
\node at (7, -5.8) {$K(Q^n_{(2)}) = \frac{4}{n}$};
\end{tikzpicture}
\caption{Graphs with constant curvature $K$: complete graphs $K_n$, cycle graphs $C_n$, hypercube graphs $Q_n$, cocktail party graphs $CP_n$,
Johnson graphs $J_{n,k}$ and demi-cube graphs $Q^n_{(2)}$.}
\end{figure}
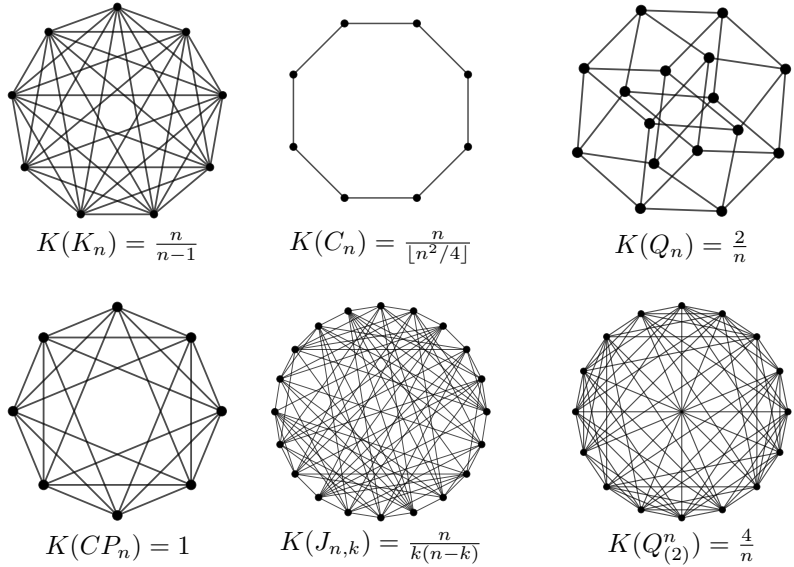
\end{center}

We start with two facts that are sometimes helpful in the computation of curvature. The first one is that positive curvature stays preserved under cartesian products.
\begin{proposition} If $G, H$ are two graphs with nonnegative curvature, then the product graph $G \times H$ also has
nonnegative curvature. If $G,H$ have constant curvature $K_1$ and $K_2$, then $G \times H$ has constant curvature $K$ where $K$ satisfies
$$ \frac{1}{K} = \frac{1}{K_1} + \frac{1}{K_2}.$$
\end{proposition}

As for the second observation, recall that a graph is vertex-transitive if, for any two vertices $v_1, v_2 \in G$, there exists an automorphism $\phi: G \rightarrow G$ such that $\phi(v_1) = v_2$. Vertex-transitive graphs always admit a constant (positive) curvature.
\begin{proposition}[Vertex-transitive graphs have constant curvature] If $G$ is vertex-transitive, then it has constant curvature $K>0$ and for any $v \in V$
$$ K =  \left( \frac{1}{n} \sum_{i=1}^{n} d(v, v_i) \right)^{-1}.$$
Any arbitrary graph $G$ with constant curvature has
$K = 1/( \frac{1}{n^2} \sum_{i,j=1}^{n} d(v_i, v_j) )^{}.$
\end{proposition}

 \textbf{Complete graph $K_n$.} The complete graph $K_n$ has constant curvature
$$ K = \frac{n}{n-1}.$$
This agrees with the Lin-Lu-Yau curvature \cite{lly} which also assigns constant curvature $n/(n-1)$ to $K_n$. It is easy to see that for any connected graph $G$ and any $D w = n \cdot \mathbf{1}$ with $\min_{i} w_i = K \geq 0$, we have
$ K \leq n/(n-1)$ with equality if and only if $G = K_n$ is the complete graph. In particular, we always have $ \min_{i} w_i \leq 2$.\\

 \textbf{Cycle graph $C_n$.} The cycle graph $C_n$ is also relatively simple: the sum over each row of $D$ is constant and for any vertex $v_i \in V$
\vspace{2pt}
$$ \sum_{j \in V} d(v_i, v_j) = \begin{cases}
\sum_{k=1}^{\left\lfloor n/2 \right\rfloor} 2k \qquad &\mbox{if}~k~\mbox{is odd}\\
\frac{n}{2} + \sum_{k=1}^{ n/2 -1} 2k \qquad &\mbox{if}~k~\mbox{is even}\\
\end{cases} \quad = \left\lfloor \frac{n^2}{4} \right\rfloor.$$
\vspace{2pt}
This implies that the cycle graph $C_n$ has constant curvature  
$$ K = \frac{n}{\left\lfloor n^2/4 \right\rfloor }  = (1+o(1)) \cdot \frac{4}{n}.$$
In contrast, both the Olivier-Ricci curvature and the Lin-Lu-Yau curvature assign curvature 0 to $C_n$ for $n \geq 6$ (see also \cite{bourne, cush1, cush2, lly2}).\\

 \textbf{Path graph $P_n$.} The trivial algebraic fact
$$ \frac{n}{n-1} (i-1) + \frac{n}{n-1} (n - i) = n$$
can be interpreted as saying that the path graph $P_n$ on $n$ vertices has curvature 0 except in the two endpoints where it has curvature $n/(n-1)$.
The Ollivier curvature of a graph is 0 on each edge while the Lin-Lu-Yau curvature is 1 on the two edges adjacent to the two endpoints and vanishes everywhere else.\\

\textbf{Hypercube graph $Q_n$.} The hypercube graph $Q_n$ with $V = \left\{0,1\right\}^n$ and edges between any two vertices with Hamming distance 1 has constant curvature
$$ K = \frac{2}{n}$$
which follows from Proposition 2 and
$$ \sum_{k=0}^{n} \binom{n}{k} k = \frac{n}{2} \cdot 2^{n}.$$
Alternatively, this would also follow from Proposition 1 and $K(Q_2) = 2$. This matches the Lin-Lu-Yau curvature \cite{lly} which also assigns constant curvature $2/n$ to $Q_n$. It is also close to the Ollivier curvature which is $2/(n+1)$ (for choice of laziness parameter $p=1/(n+1)$, we refer to Ollivier-Villani \cite{villani}).\\

 \textbf{Cocktail Party graph $CP_n$.} The cocktail party graph $CP_n$ on $2n$ vertices is defined as follows: the $2n$ vertices are split into $n$ pairs of 2 and each vertex is connected to each other vertex except the one it is paired to. For any $v \in V$, we have
 $ \frac{1}{2n} \sum_{i=1}^{2n} d(v,v_i) = 2n$
 implying that $CP_n$ has constant curvature 1 which coincides with the Ollivier-Ricci curvature.\\

 \textbf{Johnson graph $J_{n,k}$.} The Johnson graph $J_{n,k}$ is constructed as follows: the vertices are given by all $k-$element subsets of an $n-$element set. Two vertices are connected by an edge if the corresponding subsets have $k-1$ elements in common. It therefore has $\binom{n}{k}$ vertices and diameter $\min(k, n-k)$. The graphs $J_{n,k}$ and $J_{n, n-k}$ are isomorphic, we can thus assume $k \leq n/2$. The distance between two vertices $U,W$ (identified with their subsets) in $J_{n,k}$ is $d(U,W)=k - |U \cap W|$. The Johnson graph is vertex-transitive. We can thus consider the vertex $U = \left\{1,2,\dots, k\right\} \subset \left\{1,2,3,\dots, n\right\}$ and count the number of subsets $W$ with $|U \cap W| = \ell$. A moment's consideration shows that
 $$ \# \left\{W \subseteq \left\{1,2,\dots, n\right\}: |U \cap W| = k-\ell \right\} = \binom{k}{\ell} \binom{n-k}{\ell}.$$
Therefore
\begin{align*}
\sum_{W \subset \left\{1,2,\dots, n\right\}} d(U, W) &= \sum_{\ell = 0}^{k} \ell  \binom{k}{\ell} \binom{n-k}{\ell} = \frac{(n-k)k}{n} \binom{n}{k}
\end{align*}
from which, with Proposition 1, we deduce $K= n/((n-k)k)$. This again coincides with the Ollivier curvature (see \cite{rigidity} for the computation).\\

\textbf{Demi-cubes $Q^n_{(2)}$.}  $Q^n_{(2)}$ is obtained by connecting bitstrings of length $n$ if they have Hamming distance 2. This leads to two isomorphic connected components of which we pick one. This graph on $2^{n-1}$ vertices is vertex-transitive and
$$ \sum_{k=0}^{\left\lfloor n/2 \right\rfloor} k \binom{n}{2k} = \frac{n}{4} 2^{n-1}.$$
 Proposition 2 implies $K = 4/n$ which coincides with Ollivier curvature (see \cite{rigidity}). Other examples where the curvature reflects combinatorial structure of a graph in an interesting way are shown in Fig. 4.

\section{Main Results}
\subsection{An Invariant.} Several of our results will feature the quantity $\|w\|_{\ell^1}$ where $Dw = n\cdot \mathbf{1}$.  Since the linear system of equations may have multiple solutions, we start with a basic proposition for graphs admitting nonnegative curvature.
\begin{proposition}[Invariance of total curvature] Let $G$ be a connected graph and suppose $Dw_1 = n\cdot \mathbf{1} = D w_2$ for two vectors $w_1, w_2 \in \mathbb{R}_{\geq 0}^{n}$. Then $\|w_1\|_{\ell^1} = \|w_2\|_{\ell^1}$.
\end{proposition}
The quantity $\|w\|_{\ell^1}$, the sum over all curvatures, plays a role in many of our results.  Proposition 3 guarantees that the results do not depend on which solution of $Dw = n\cdot \mathbf{1}$ (should multiple exist) one chooses. There is an interesting subtlety to Proposition 3: it appears as if it should have a simple proof via linear algebra. After all, if two different $w_1, w_2 \in \mathbb{R}_{\geq 0}^{n}$ with  $Dw_1 = n\cdot \mathbf{1} = D w_2$ exist, then their difference is in the nullspace $w_1 - w_2 \in \ker(D)$ and if the nullspace of $D$ is orthogonal to the constant vector $(1,1,\dots, 1)$, we have the desired result. This, however, is not always the case: there are graphs $G$ (see Fig. 5) for which the nullspace of the associated distance matrix $D$ may not have this property. However, in those cases $D w = n \cdot \mathbf{1}$ will not have a solution $w \in \mathbb{R}^n_{\geq 0}$. Conversely, Proposition 3 guarantees that if $D w = n \cdot \mathbf{1}$ has a solution in $w \in \mathbb{R}^n_{\geq 0}$, then the nullspace will either be empty or orthogonal to the constant vector.

\subsection{Discrete Bonnet-Myers theorem} The classical Bonnet-Myers theorem \cite{myers} states that if $(M,g)$ is a complete, connected $n-$dimensional manifold with Ricci curvature bounded below by $1/r^2$, then $\diam(M) \leq \pi r$: a manifold with uniformly positive curvature cannot be too large. Cheng \cite{cheng} later proved that equality can only happen in the case of constant sectional curvature. The same type of result holds for our notion of curvature on graphs.

\begin{theorem}[Discrete Bonnet-Myers] Let $G$ be a connected graph. 
If $Dw = n \cdot \mathbf{1}$ has curvature bounded from below by $K = \min_i w_i  \geq 0$, then
$$ \diam(G) \leq \frac{2n}{\|w\|_{\ell^1}} \leq \frac{2}{K}.$$
If $\diam(G) \cdot K = 2$, then $G$ has constant curvature. 
\end{theorem}

This result is sharp: examples are given by even cycles $C_{2n}$, the hypercube graphs $Q_n$ or the Johnson graph $J_{2n,n}$. Indeed, there are many examples for which $\diam(G) = 2/K$, some are shown in Fig 3.  Theorem 1 matches the discrete Bonnet-Myers Theorem that has been established for Ollivier-Ricci curvature \cite{bourne, lly, ollivier}.

\begin{center}
\begin{figure}[h!]
\begin{tikzpicture}
\node at (0,0) {\includegraphics[width=0.22\textwidth]{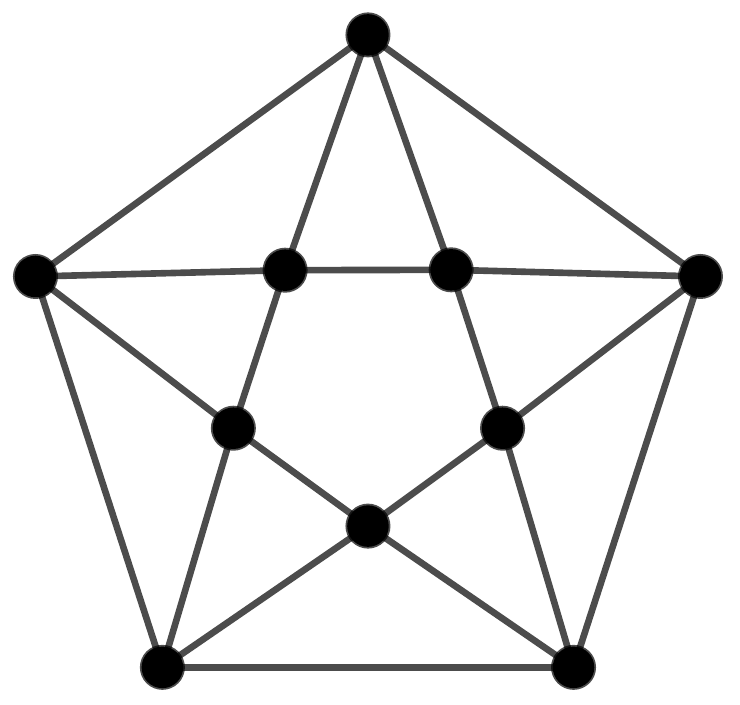}};
\node at (3.5,0) {\includegraphics[width=0.22\textwidth]{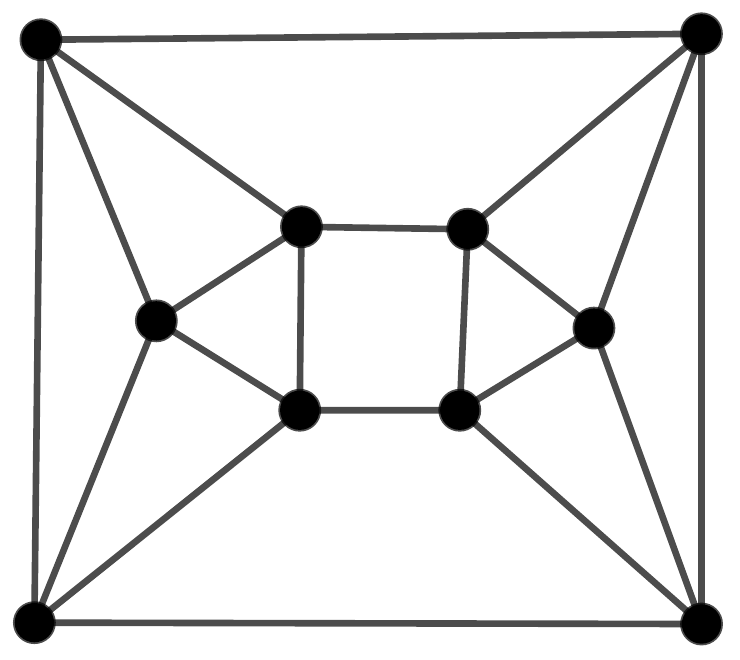}};
\node at (7,0) {\includegraphics[width=0.22\textwidth]{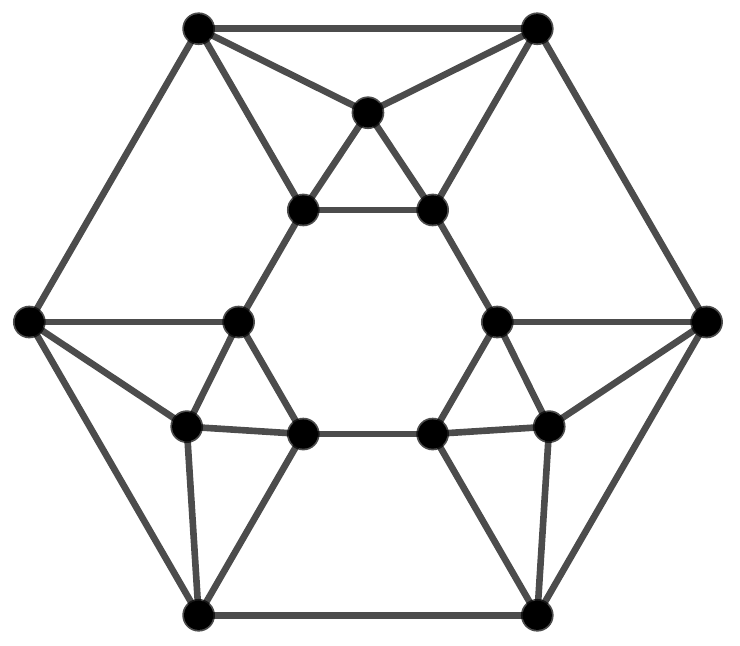}};
\node at (0,-3) {\includegraphics[width=0.22\textwidth]{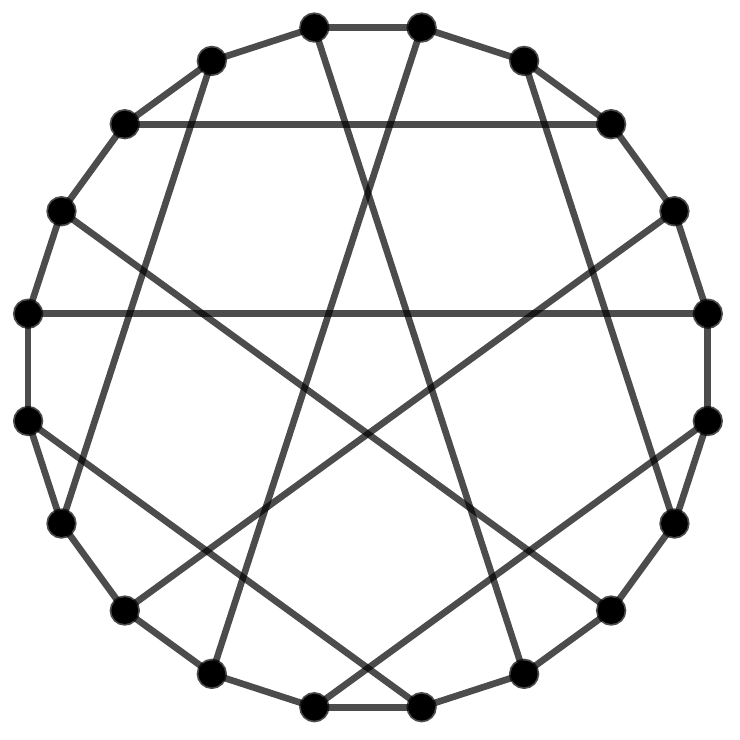}};
\node at (3.5,-3) {\includegraphics[width=0.22\textwidth]{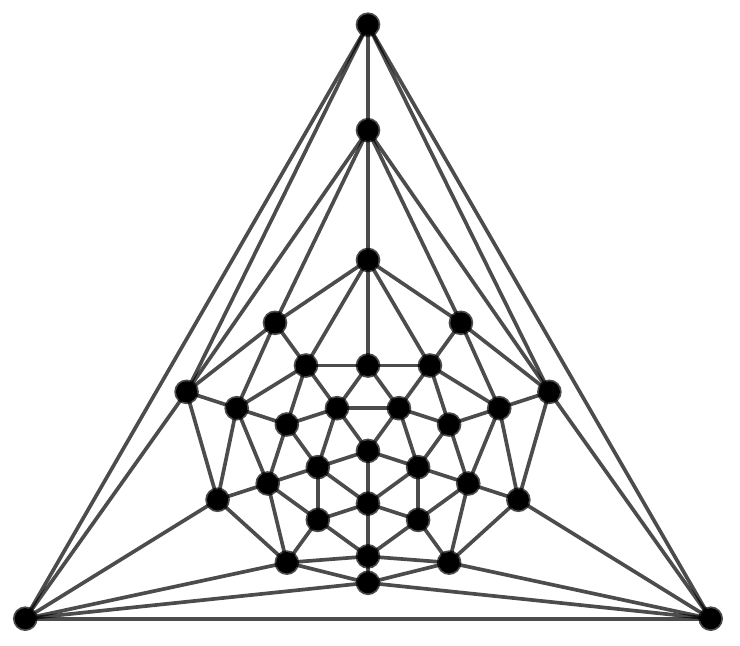}};
\node at (7,-3) {\includegraphics[width=0.22\textwidth]{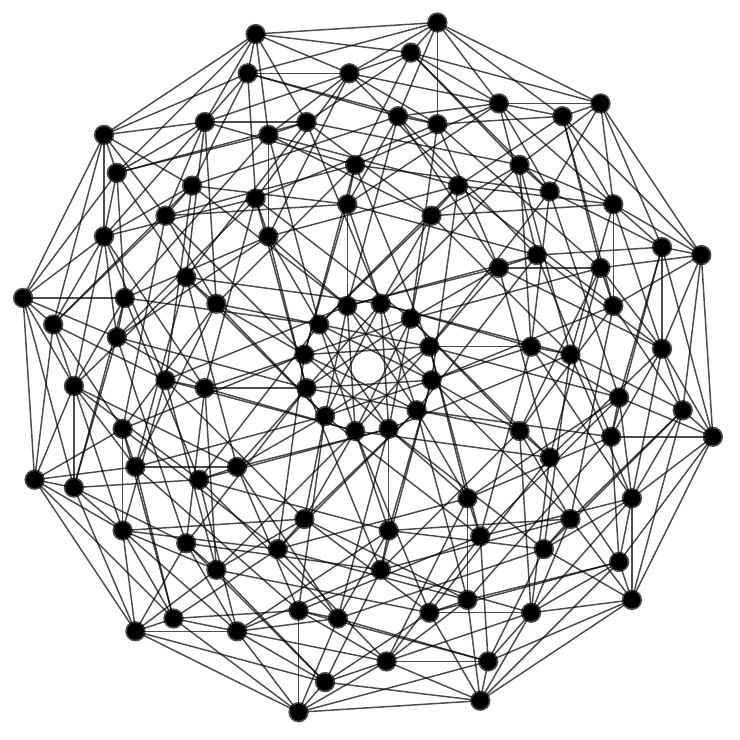}};
\end{tikzpicture}
\caption{Examples of graphs for which Theorem 1 is sharp and $K \cdot \diam(G) = 2$ (all have constant curvature).}
\end{figure}
\end{center}

We also obtain a result in the other direction: if a graph is nonnegatively curved and the diameter is small, then this forces the total curvature to be large.
\begin{theorem}[Reverse Bonnet-Myers] Let $G$ be a graph with nonnegative curvature, i.e. there exists $Dw = n \cdot \mathbf{1}$ with $\min_i w_i  \geq 0$. Then
$$ \|w\|_{\ell^1}  \geq  \frac{n^2}{n-1}\frac{1}{\diam(G)}$$
with equality if and only if $G=K_n$.
\end{theorem}
We note a particular consequence of Theorem 1 and Theorem 2: if $G$ has constant curvature $K>0$, then the curvature is inversely proportional to the diameter and
$$\frac{1}{K} \leq \diam(G) \leq \frac{2}{K}.$$

\subsection{Discrete Lichnerowicz Theorem} 
The classic Lichnerowicz Theorem \cite{lich} states that if an $n-$dimensional manifold has Ricci curvature bounded below by $K > 0$, then the first eigenvalue of the Laplacian satisfies $\lambda_1 \geq n/(n-1) K$. A natural analogue of the first eigenfunction on a graph is the smallest nontrivial eigenvalue of the Laplacian matrix $L = D - A$ which can be defined as the minimum of a quadratic form over all functions with mean value 0 or, equivalently,
$$ \lambda_1(G)  = \inf_{f:V \rightarrow \mathbb{R} \atop \sum_{v \in V} f(v) = 0} \frac{ \sum_{(u,v) \in E}~ (f(u) - f(v))^2 }{\sum_{v \in V} f(v)^2}.$$
\begin{theorem}[Discrete Lichnerowicz] Let $G$ be a connected graph. If $Dw = n \cdot \mathbf{1}$ has curvature bounded from below $K = \min_i w_i  > 0$, then
$$ \lambda_1 \geq \frac{ \|w\|_{\ell^1}}{2 n^2} \geq  \frac{K}{2n}.$$
\end{theorem}
The result is sharp up to constants: if we take the cycle graph $C_n$, then
$$ \lambda_1 = 4 \sin\left( \frac{\pi}{n}\right)^2 \sim \frac{4 \pi^2}{n^2} \qquad \mbox{while} \qquad \frac{K}{2n} =  \frac{1}{2 \left\lfloor n^2/4 \right\rfloor }  = (1+o(1)) \cdot \frac{2}{n^2}.$$
We note that the Lichnerowicz scaling for other notions of curvature (say, Ollivier curvature or Lin-Lu-Yau curvature) ends up being different since, in that case, $\lambda_1 \geq K$. As the cycle graph shows, this is clearly not possible here. However, the cycle graph does play a somewhat distinguished role: a result of Lin-Lu-Yau \cite{lly2} (see also \cite{cush1, cush2}) implies that if a finite graph with girth at least 5 and $n > 20$ vertices has vanishing Lin-Lu-Yau curvature, then $G=C_n$.  It follows from Theorem 2 that if a graph on $n$ vertices is nonnegatively curved $\min_i w_i \geq 0$, then $\max_i w_i \geq 1/\diam(G) \geq 1/n$. It seems conceivable that the sharp bound for this estimate might actually be $\max_i w_i \geq n/\left\lfloor n^2/4\right\rfloor$ with equality if and only if $G=C_n$.

\subsection{Total Curvature Minimax.} Let $G$ be nonnegatively curved, i.e. assume $Dw = n \cdot \mathbf{1}$ has a solution with $\min_i w_i \geq 0$. Then $G$ has the following interesting balancing property: for any (weighted) collection of vertices, there always exists another vertex $a \in V$ such that the average distance between $a$ and our collection of vertices is not too large. Moreover, there also exists a vertex $b \in V$ such that the average distance between the collection and $b$ is not too small. This is the strongest Theorem in this paper: it is then used to prove all previous results.
 
\begin{theorem}[Minimax Theorem]
Let $G$ be nonnegatively curved with total curvature $\|w\|_{\ell^1}$.  Then, for \emph{any} probability measure $\nu$ on $V$, there are $a, b \in V$ with
$$ \min_{a \in V}  \sum_{v \in V}^{} d(a, v) \nu(v) \leq \frac{n}{\|w\|_{\ell^1}} \leq \max_{b \in V}  \sum_{v \in V}^{} d(b, v) \nu(v).$$
\end{theorem}
We emphasize that $n/\|w\|_{\ell^1}$ is an invariant of the graph and completely independent of the measure $\nu$.
Note also that, in particular, if $G$ has constant curvature $K$, then the result implies the existence of $a,b \in V$ with
$$  \sum_{v \in V}^{} d(a, v) \nu(v) \leq \frac{1}{K} \leq  \sum_{v \in V}^{} d(b, v) \nu(v).$$
$n/\|w\|_{\ell^1}$ is the unique number with this property: if $\nu = w/\|w\|_{\ell^1}$, both inequalities are sharp. This is a consequence of the von Neumann Minimax Theorem.

To the best of our knowledge, this kind of property has not been considered for any of the other notions of curvature. Given the delicate nature of the statement, one would perhaps not expect it to hold in general but it could be interesting to understand whether approximate versions for other types of curvature hold true.
In the context of connected, bounded metric spaces, such results date back to a 1964 paper of Gross \cite{gross} who showed that there exists a number $\alpha > 0$ (the `rendezvous number of the metric space') such that for any (weighted) collection of points in the space there always exists another point at average distance exactly $\alpha$. We refer to the survey of Cleary \& Morris \cite{cleary}. Gross' theorem, which originally appeared in \textit{Advances in Game Theory}, makes use of
a 1952 result of Glicksberg \cite{glicksberg} which `implies the minimax theorem for continuous games with continuous payoff as well as the existence of Nash equilibrium points'. In light of this, it is perhaps less surprising that the von Neumann Minimax theorem \cite{john} would appear. The existence of such an $\alpha$ for finite metric spaces was also shown by Thomassen \cite{thom}.

\subsection{Inverting the Linear System.}

The equation  $D w = n \cdot \mathbf{1}$ need not always have a solution, however, this seems to be exceedingly rare. 
 In a search of all 9059 graphs with $2 \leq n \leq 500$ vertices that are implemented in Mathematica, there are five examples where the linear system does not have a solution (listed in Table 1).

\begin{center}
\begin{table}[h!]
\begin{tabular}{l | c | c | c | c }
Graph & $\#V$  & $\#E$  &     $  D^{\dagger} (n \cdot \mathbf{1}) \subseteq$   &   $  D(D^{\dagger} (n \cdot \mathbf{1})) \subseteq$ \\
\hline
\vspace{3pt}
$K_{1,1,1,4}$    & 7 & 15 &      $ \left[0.65, 0.99 \right]$    &       $ \left[5.25, 7.875 \right]$    \\
\vspace{3pt}
$K_{1,1,1,1,3}$    & 7 & 18 &      $  \left[0.85, 1.15\right]$    &       $ \left[6, 8 \right]$    \\
\vspace{3pt}
Quartic$-(11,18)$    & 11 & 22 &      $ \left[0.16, 1.05 \right]$    &       $ \left[10.32, 11.40 \right]$    \\
\vspace{3pt}
Cubic$-(14,52)$   & 14 & 21 &      $ \left[-1.09, 2.22  \right]$    &       $ \left[13.02, 14.97 \right]$    \\
\vspace{3pt}
Knight$-(7,7)$   & 49 & 120 &      $ \left[-10.93, 2.75  \right]$    &       $ \left[46.42, 52.22 \right]$    \\
\end{tabular}
\caption{Five exceptional graphs: $Dv = n \cdot \textbf{1}$ does not have a solution.}
\end{table}
\end{center}
\vspace{-30pt}
In such cases, the Moore-Penrose pseudo-inverse $w =D^{\dagger}(n \cdot \mathbf{1})$ is a natural replacement: recall that the pseudo-inverse is the vector $z$ minimizing $\| Dz - n \cdot \mathbf{1} \|_{\ell^2}$. If there is more than one such vector, then it is defined as the one with smallest $\ell^2-$norm which is then uniquely determined. For the five exceptional cases the linear system can \textit{almost} be solved: the vector $D(D^{\dagger}(n \cdot \mathbf{1}))$ is nearly constant. 
\begin{theorem}[Discrete Bonnet-Myers and Lichnerowicz II] Let $G$ be a connected graph and let $w \in \mathbb{R}^n_{>0}$ be arbitrary. Then, for $K = \min_i w_i > 0$, we have
$$ \diam(G) \leq \frac{\|Dw\|_{\ell^{\infty}}}{n} \frac{8}{K} $$
and
$$ \lambda_1 \geq \frac{1}{8 \|Dw\|_{\ell^{\infty}}}  K.$$
\end{theorem}

In the case where $Dw = n \cdot \mathbf{1}$ has a solution, we have $\|Dw\|_{\ell^{\infty}} = n$ and recover the bounds $\diam(G) \leq 8/K$ and $\lambda_1 \geq K/(8n)$ which are optimal up to constants. Theorem 5 can be applied to the first three of the five exceptional graphs in Table 1.
At this point one could wonder how many exceptional graphs there are: for which graphs does $Dw = n \cdot \mathbf{1}$ not have a solution? While such graphs exist, they seem to be rather rare and their relative proportion seems to decrease as the number of vertices increases: taking thousands of Erd\H{o}s-Renyi graphs, one can find examples with $n=10$ vertices but, using random sampling, the proportion of such examples seems to rapidly decrease as $n$ increases. Moreover, in all the examples we found that for $w = D^{\dagger}(n \cdot \mathbf{1})$ all the entries of $Dw$ are approximately constant in the sense that they all are in the $[0.7n, 1.3n]$ range. This seems like an interesting question in its own right: is there a reason why $Dw = n \cdot \mathbf{1}$ seems to almost always have a solution?  Is there a reason why exceptional graphs are rare? Is there exists a constant $0 < c < 1$ such that for all graphs $  D(D^{\dagger} (n \cdot \mathbf{1})) \subset [cn, c^{-1} n]$?

 \begin{center}
\begin{figure}[h!]
\begin{tikzpicture}[scale=1]
\node at (0,3.5) {\includegraphics[width=0.24\textwidth]{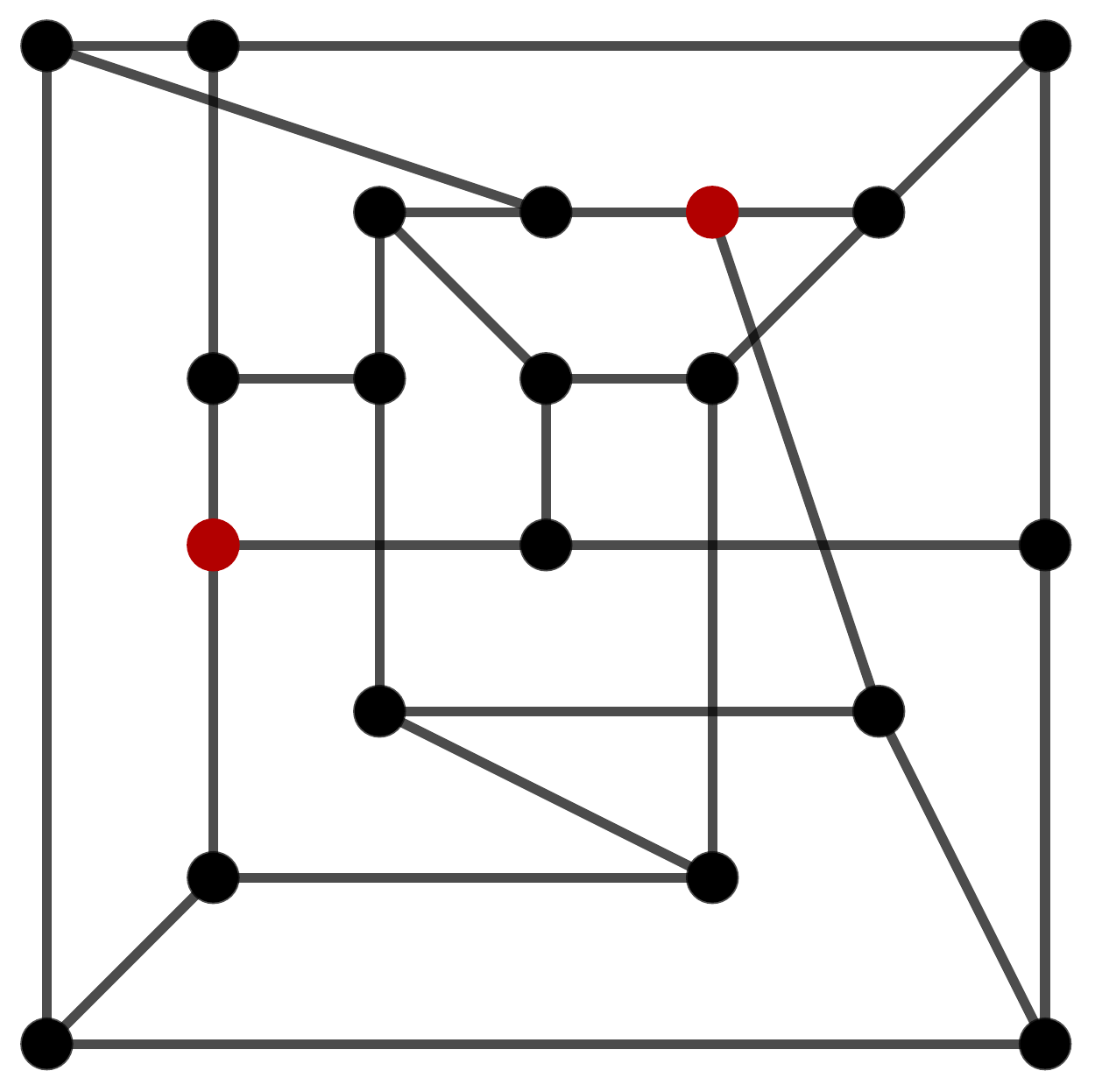}};
\node at (3.5,3.5) {\includegraphics[width=0.24\textwidth]{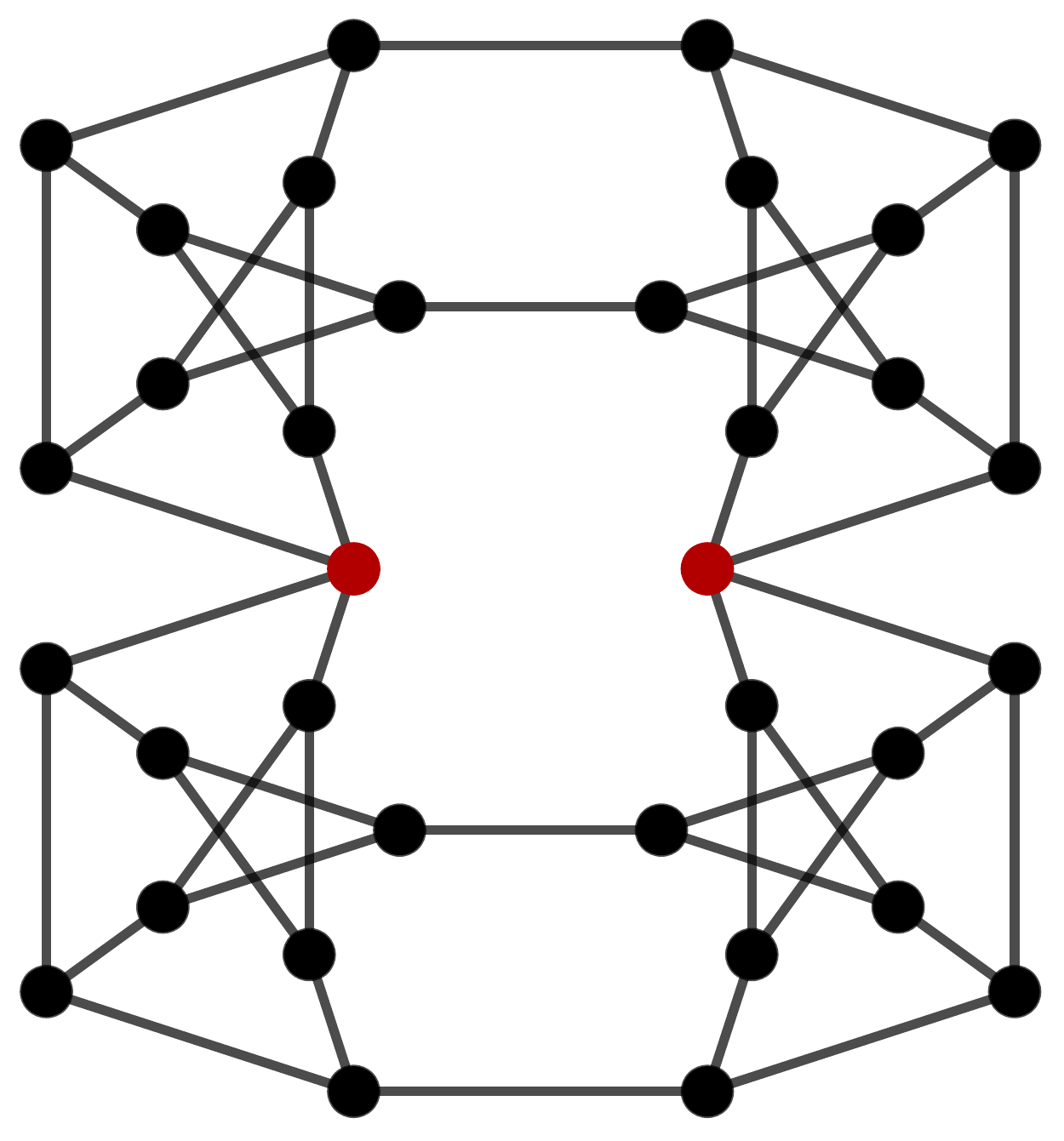}};
\node at (7,3.5) {\includegraphics[width=0.25\textwidth]{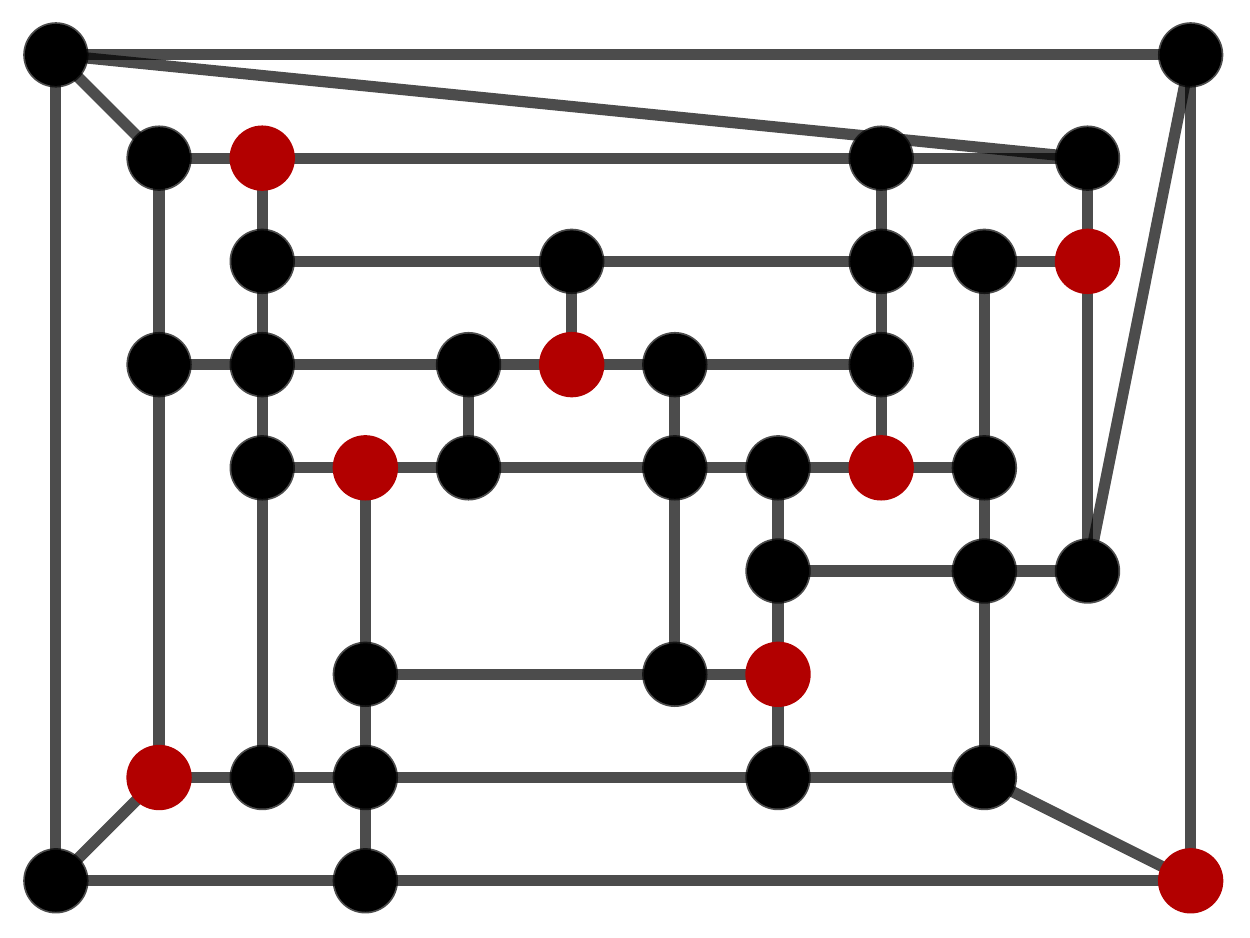}};
\node at (0,0) {\includegraphics[width=0.3\textwidth]{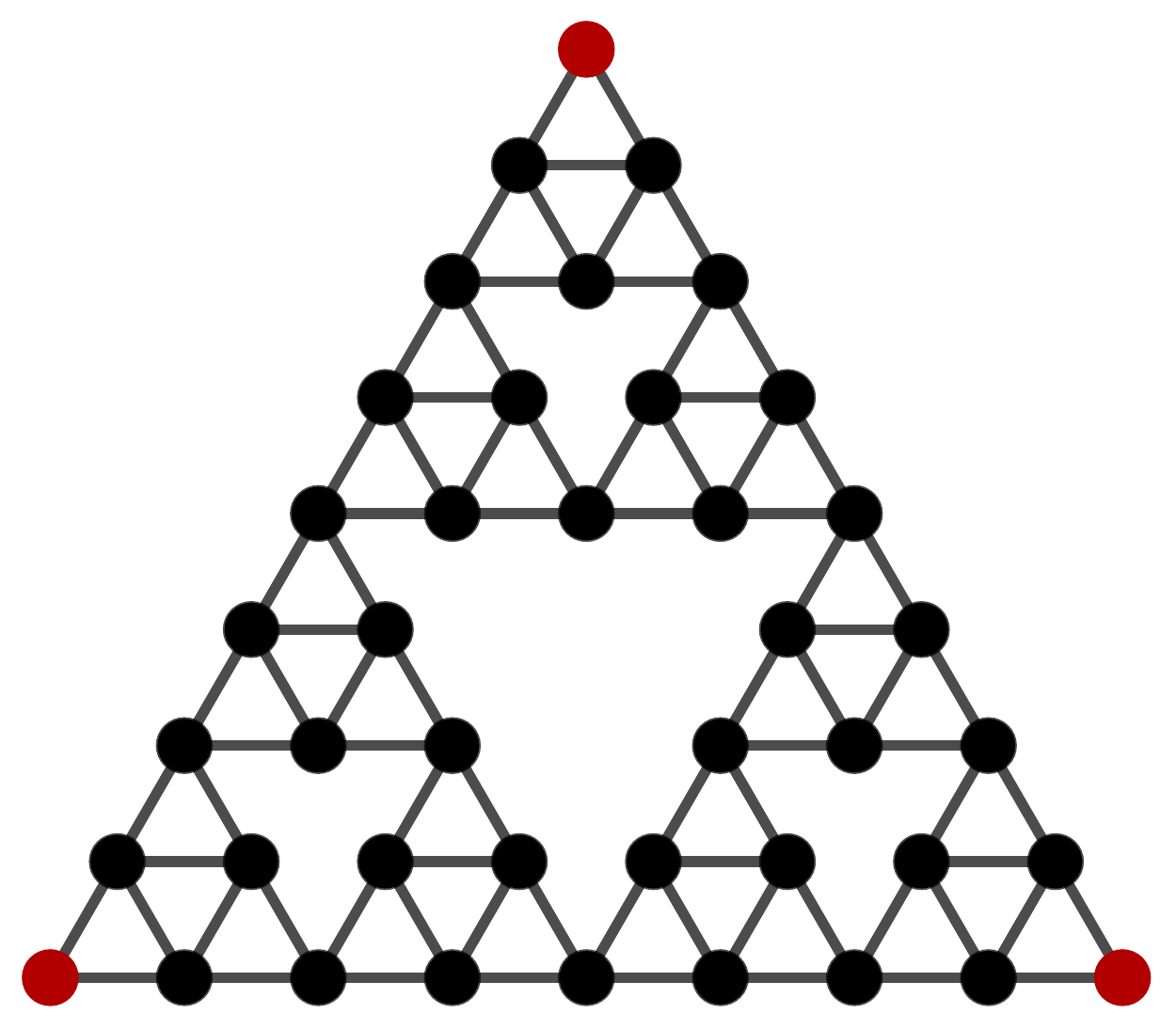}};
\node at (5,0) {\includegraphics[width=0.3\textwidth]{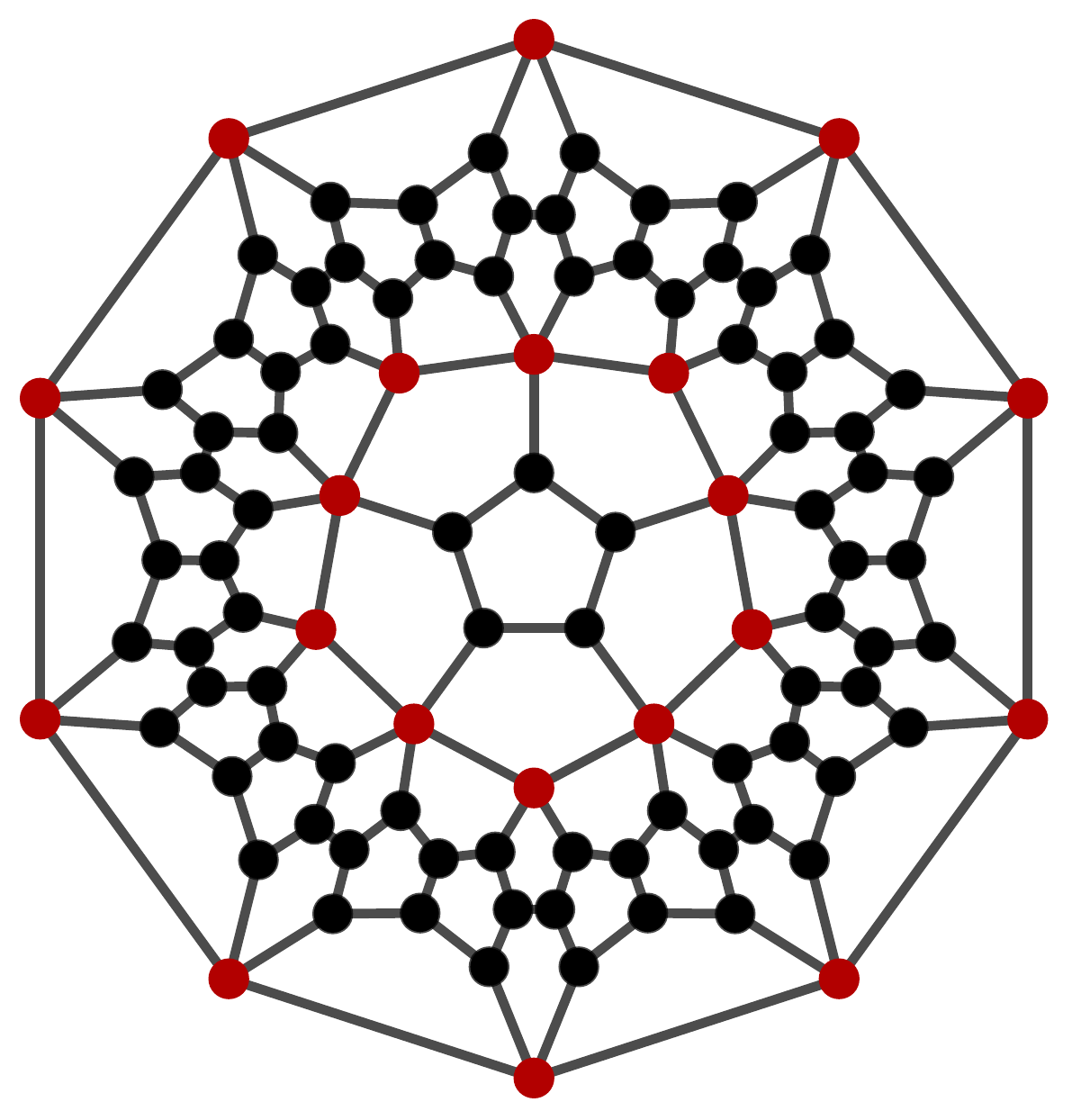}};
\end{tikzpicture}
\vspace{-10pt}
\caption{Examples of nonnegatively curved graphs with positive curvature concentrated in a few isolated vertices (red) while most vertices have curvature 0 (black). Top: CrossingNumberGraph6C, Thomassen 34 and the Pentagonal Icositetrahedral Graph, Bottom:
Sierpinski and Thomassen 105.}
\end{figure}
\end{center}

We have since investigated the phenomenon in subsequent work \cite{stein}. There, the following sufficient criterion for the existence of a solution $D w = \mathbf{1}$ was given.
\begin{prop}[\cite{stein}] Suppose $D \in \mathbb{R}_{\geq 0}^{n \times n}$ has eigenvalues $\lambda_1 > 0 \geq \lambda_2 \geq \dots \geq \lambda_n$ and eigenvector $Av = \lambda_1 v$. If 
$$  1 - \left\langle v, \frac{\mathbf{1}}{\sqrt{n}}\right\rangle^2 < \frac{|\lambda_2|}{\lambda_1 - \lambda_2},$$
then $Dx = \mathbf{1}$ has a solution.
\end{prop}
The applicability of this Proposition depends on two things: how often a graph distance matrix $D$ has such a spectral structure and how often the inequality is satisfied. It seems that the spectral structure is fairly common but not so
common as to explain the overall solvability of $D x = \mathbf{1}$ (it is only a sufficient criterion). However, in what is somewhat remarkable, the inequality tends to be satisfied quite frequently. This leads to a related phenomenon which we
describe as follows.

\begin{quote}
\textbf{Phenomenon.} Let $G=(V,E)$ be a connected, bounded graph, let $D \in \mathbb{R}^{n \times n}$ be its distance matrix and let $v \in \mathbb{R}^{n}_{\geq 0}$ denote the eigenvector corresponding to the largest eigenvalue of $D$ (which, by Perron-Frobenius, can be taken to have non-negative entries). Then this eigenvector is `nearly constant' in the sense that 
$$ c_G = \frac{ \left\langle v, \mathbf{1} \right\rangle}{   \| v\|_{\ell^2}  \cdot \|\mathbf{1} \|_{\ell^2}}$$
tends to be very close to 1. 
\end{quote}

By Cauchy-Schwarz, we have $c_G \leq 1$.  The main result of \cite{stein} ensures that, for general distance matrices in metric spaces, $c_G \geq 1/\sqrt{2}$. It seems that for most graphs, the constant is much close to 1. Indeed, it is not easy to find examples where $c_G \leq 0.95$, we refer to \cite{stein} for details. In summary, it seems that the solvability of $Dw = \mathbf{1}$ is connected to a number of interesting phenomena which, at present, are poorly understood and suggest interesting avenues for further research.

\section{Proofs}
We first establish Proposition 1 (in \S 3.1) and Proposition 2 (in \S 3.2). These two arguments are elementary and self-contained. \S 3.3 is concerned with the proof of Theorem 4 (the Minimax Theorem) which is perhaps the most substantial result in the paper and will be used to prove most subsequent results. \S 3.4 then uses Theorem 4 to prove Theorem 1 (Bonnet-Myers) and Theorem 2 (Reverse Bonnet-Myers). \S 3.5 establishes Proposition 3, the invariance of total curvature for positively-curved graphs, this argument will also be a consequence of Theorem 4. \S 3.6 proves the Lichnerowicz inequality (Theorem 3) which follows quickly from Theorem 1 and a standard spectral estimate (whose short proof is included for the convenience of the reader). Finally, \S 3.7 establishes, by a completely independent combinatorial argument that is unrelated to any prior arguments, Theorem 5.

\subsection{Proof of Proposition 1}
\begin{proof}
Let $G$ be a graph with curvature bounded from below by $K_1 \geq 0$, meaning there exists $w_1 \in \mathbb{R}^{n_1}_{\geq 0}$ with $D w_1 = n_1$ and $K_1 = \min_i (w_1)_i$, and let $H$ be another graph with curvature bounded from below by $K_2 \geq 0$, meaning there exists $w_2 \in \mathbb{R}^{n_2}_{\geq 0}$ with $D w_2 = n_2$ and $K_2 = \min_i (w_2)_i$. Our goal is to show that $G \times H$ has curvature bounded from below by 0.
We will consider the product measure $\mu = w_1 \times w_2$. Then, for any $(g_1, h_1) \in V(G \times H)$, the expression
$$ X =  \sum_{(g_2,h_2) \in V(G \times H)} d((g_1, h_1), (g_2,h_2)) \mu((g_2, h_2))$$
can be simplified to
\begin{align*}
X&=  \sum_{(g_2,h_2) \in V(G \times H)} d((g_1, h_1), (g_2,h_2)) w_1(g_2) w_2(h_2) \\
 &=  \sum_{(g_2,h_2) \in V(G \times H)} \left(d(g_1, g_2) + d(h_1, h_2)\right) w_1(g_2) w_2(h_2). 
 \end{align*}
 At this point, the sum factors into two sums and
 \begin{align*}
X &=\sum_{(g_2,h_2) \in V(G \times H)} d(g_1, g_2) w_1(g_2) w_2(h_2) \\
 &+ \sum_{(g_2,h_2) \in V(G \times H)} d(h_1, h_2) w_1(g_2) w_2(h_2) = \|w_2\|_{\ell^1} \cdot n_1 + \|w_1\|_{\ell^1}\cdot n_2
 \end{align*}

 which is a positive constant independently of $(g_1, h_1) \in G \times H$. This means that the rescaled measure
 \begin{align*}
  w &= \frac{n_1 n_2}{\|w_2\|_{\ell^1} \cdot n_1 + \|w_1\|_{\ell^1}\cdot n_2} (w_1 \times w_2)\\
  &=  \left( \frac{\|w_1\|_{\ell^1}}{n_1} + \frac{ \|w_2\|_{\ell^1}}{n_2} \right)^{-1} (w_1 \times w_2)
  \end{align*}
satisfies $D_{G \times H} w = n_1 n_2 = |V(G \times H)|$ and thus is an admissible nonnegative curvature on $G \times H$. Moreover,
 $$ \min w =  \left( \frac{\|w_1\|_{\ell^1}}{n_1} + \frac{ \|w_2\|_{\ell^1}}{n_2} \right)^{-1} K_1 K_2 \geq 0.$$
If $G$ and $H$ have constant curvature $K_1$ and $K_2$ then
  $$ K =  \min w =  \left( K_1 + K_2\right)^{-1} K_1 K_2$$
  and thus
  $$ \frac{1}{K} = \frac{1}{K_1} + \frac{1}{K_2}.$$
\end{proof}

We note the following immediate consequence.

\begin{corollary}
If $G$ has constant curvature $K$, then $G^n = G \times \dots \times G$ has constant curvature $K/n$.
\end{corollary}

This statement has a direct analogue for Ollivier-Ricci curvature and Lin-Lu-Yau curvature.  Indeed, the Corollary
under the additional assumption of $G$ being regular, is true verbatim for Lin-Lu-Yau curvature, see \cite[Corollary 3.2]{lly}.

\subsection{Proof of Proposition 2}
\begin{proof}
If $G$ is vertex-transitive, then for each $u \in V$ and each $k \in \mathbb{N}$, the size of
$$ \# \left\{v \in V: d(u,v) = k\right\} \qquad \mbox{only depends on}~k~\mbox{and not on}~u.$$
This means that the rows of the distance matrix $D$ are permutations of each other and, in particular, that for each vertex $v_i \in V$ the row sum
$$ \sum_{j=1}^{n} d(v_i, v_j)  = R \qquad \mbox{is independent of}~i.$$
This shows that the graph admits a constant positive curvature with $K= n/R$. Assume now that $G$ admits constant curvature $K$. Then, for each $1 \leq i \leq n$, 
$$ \sum_{j=1}^{n} d(v_i, v_j) K = n$$
and the result follows by summing over $i$.
\end{proof}

\subsection{Proof of Theorem 4}
 Theorem 4 uses a specific implication of the von Neumann Minimax Theorem which reads as follows.

\begin{thm}[von Neumann \cite{john}] Let $A \in \mathbb{R}^{n \times n}$ by a symmetric matrix. There exists a unique $\alpha \in \mathbb{R}$ such that for all $(x_1, \dots, x_n) \in \mathbb{R}^n_{\geq 0}$ satisfying $x_1 + \dots + x_n = 1$
$$ \min_{1\leq i \leq n} (Ax)_i \leq \alpha \leq \max_{1\leq i \leq n} (Ax)_i.$$ 
\end{thm}
Since the statement deviates a little from the way the Minimax theorem is usually phrased, we quickly deduce it from the more canonical formulation.
\begin{proof} The way the Minimax theorem is typically phrased (for quadratic matrices) is as follows: given an arbitrary matrix $A \in \mathbb{R}^{n \times n}$, we consider the space of mixed strategies for both players
$$ X =  \left\{z \in \mathbb{R}^n: \forall ~1 \leq i \leq n: z_i \geq 0 \quad \mbox{and} \quad \sum_{i=1}^{n} z_i = 1\right\} = Y,$$
where $X$ are the strategies that can be played by Player 1 and $Y$ are the strategies that can be played by Player 2. The pay-off of any given pair of strategies $(x,y) \in X \times Y$ is $x^T A y = \left\langle x, Ay \right\rangle$. The goal of Player 1 is to maximize the pay-off while the goal of Player 2 is to minimize the pay-off.
The Minimax Theorem then states that the game has a value $\alpha \in \mathbb{R}$ which means that
$$ \max_{x \in X} \min_{y \in Y} \left\langle x, Ay \right\rangle = \alpha =  \min_{y \in Y} \max_{x \in X} \left\langle x, Ay \right\rangle.$$
The first equation implies that there exists a strategy $x^* \in X$ such that Player 1 can always guarantee payoff at least $\alpha$ independently
of what Player 2 is doing. The second equation implies the existence of a strategy $y^* \in Y$ such that Player 2 can always guarantee a pay-off of at most $\alpha$ independently of what Player 1 is doing. We will now consider additionally that $A$ is symmetric.\\

Note that, for any given action by Player 2, a fixed $y \in Y$, it is clear how Player 1 would react: they would select the largest pay-off (which may or may not be unique). This means that, for any fixed $y \in Y$,
$$  \max_{x \in X} \left\langle x, Ay \right\rangle = \max_{1 \leq i \leq n} (Ay)_i,$$
where $(Ay)_i$ denotes the $i-$th entry of the vector and therefore
$$ \min_{y \in Y} \max_{x \in X} \left\langle x, Ay \right\rangle = \min_{y \in Y} \max_{1 \leq i \leq n} (Ay)_i.$$
Using the symmetry of the matrix, we can use the same logic to write
$$ \max_{x \in X} \min_{y \in Y} \left\langle x, Ay \right\rangle =  \max_{x \in X} \min_{y \in Y} \left\langle Ax, y \right\rangle = \max_{x \in X} \min_{1 \leq i \leq n} (Ax)_i.$$
Altogether, we arrive at
$$  \max_{x \in X} \min_{1 \leq i \leq n} (Ax)_i = \alpha =  \min_{x \in X} \max_{1 \leq i \leq n} (Ax)_i.$$
It now follows that for any arbitrary linear combination of the rows $z \in X$
$$ \min_{1 \leq i \leq n} (Az)_i \leq \max_{x \in X} \min_{1 \leq i \leq n} (Ax)_i = \alpha = \min_{x \in X} \max_{1 \leq i \leq n} (Ax)_i  \leq \max_{1 \leq i \leq n} (Az)_i.$$
\end{proof}

\begin{proof}[Proof of Theorem 4] We use the von Neumann Minimax Theorem when applied to the distance matrix $D$ of a graph $G$. Let $\nu$ be an arbitrary probability measure on the vertices. Then the function $f:V \rightarrow \mathbb{R}$ given by
$$ f(u) = \sum_{v \in V} d(u,v) \nu(v)$$
can be written as a vector provided we interpret the measure as a vector $\nu \in \mathbb{R}^n$ 
$$ f = D \nu.$$
The Minimax Theorem now implies the existence of a unique number $\alpha \in \mathbb{R}$ independent of the measure $\nu$ such that
$$ \min_{1 \leq i \leq n} (D \nu)_i \leq \alpha \leq  \max_{1 \leq i \leq n} (D \nu)_i.$$
Rewriting this in terms of distances, this merely says
$$ \min_{a \in V} \sum_{v \in V}^{} d(a, v) \nu(v) \leq \alpha \leq \max_{b \in V}  \sum_{v \in V}^{} d(b, v) \nu(v).$$
If $Dw = n \cdot \mathbf{1}$ for $w \in \mathbb{R}^{n}_{\geq 0}$, then we consider the probability measure 
$$\nu = \frac{w}{\|w\|_{\ell^1}}$$ 
and find that for each $u \in V$
$$ \sum_{v \in V}^{} d(u, v) \nu(v) = \frac{1}{\|w\|_{\ell^1}} \sum_{v \in V}^{} d(u, v)  w_v = \frac{(D w)_{u}}{\|w\|_{\ell^1}} = \frac{n}{\|w\|_{\ell^1}}.$$
Since this is true for every individual vertex, it is also true for the maximum and the minimum and we may deduce that the unique number has to satisfy
$$ \alpha =  \frac{n}{\|w\|_{\ell^1}}.$$
\end{proof}

\subsection{Proof of Theorem 1 and Theorem 2}
\begin{proof} The Bonnet-Myers and reverse Bonnet-Myers theorem now follow quickly from Theorem 4.
Suppose $Dw = n \cdot \mathbf{1}$ and $\min_i w_i = K \geq 0$. Take two vertices $a,b \in V$ at maximal distance $d(a,b) = \diam(G)$. We apply Theorem 4 to the set of vertices $\left\{a,b\right\}$ with $\nu(a) = \nu(b) = 1/2$ and conclude that there exists $c \in V$ with
$$\sum_{v \in V}^{} d(c, v) \nu(v)= \frac{1}{2} \left( d(a,c) + d(b,c) \right) \leq \frac{n}{\|w\|_{\ell^1}}.$$
Using the triangle inequality, we have that
$$ \frac{1}{2} \left( d(a,c) + d(b,c) \right) \geq \frac{1}{2} d(a,b) = \frac{\diam(G)}{2}.$$
This implies Theorem 1. If we pick $\nu$ to be the uniform probability measure, then Theorem 4 implies that the existence of a vertex $b$ such that
$$  \frac{n}{\|w\|_{\ell^1}} \leq \sum_{v \in V}^{} d(b, v) \nu(v)= \frac{1}{n} \sum_{v \in V}^{} d(b, v).$$
Considering $d(b,b) = 0$, we have
\begin{align*}
 \frac{1}{n} \sum_{v \in V}^{} d(b, v) =  \frac{1}{n} \sum_{v \in V \atop v \neq b}^{} d(b, v) \leq  \frac{1}{n} \sum_{v \in V \atop v \neq b}^{} \diam(G) = \frac{n-1}{n} \diam(G).
 \end{align*}
 This implies the desired result. Equality can only happen if, for all $v \neq b$, we have $d(b,v) = \diam(G)$ which then implies that $G = K_n$.
\end{proof}

\subsection{Proof of Proposition 3}

We note that one natural property that would imply Proposition 3 would be if the nullspace of a distance matrix $D$ only
contained vectors whose entries add up to 0. This, however, is not always the case (even though exceptions seem to be exceedingly rare). Two counterexamples are shown in Fig. 5. 
Proposition 3 implies that these graphs do not admit positive curvature. 
 \begin{center}
\begin{figure}[h!]
\begin{tikzpicture}[scale=1]
\node at (0,3.5) {\includegraphics[width=0.2\textwidth]{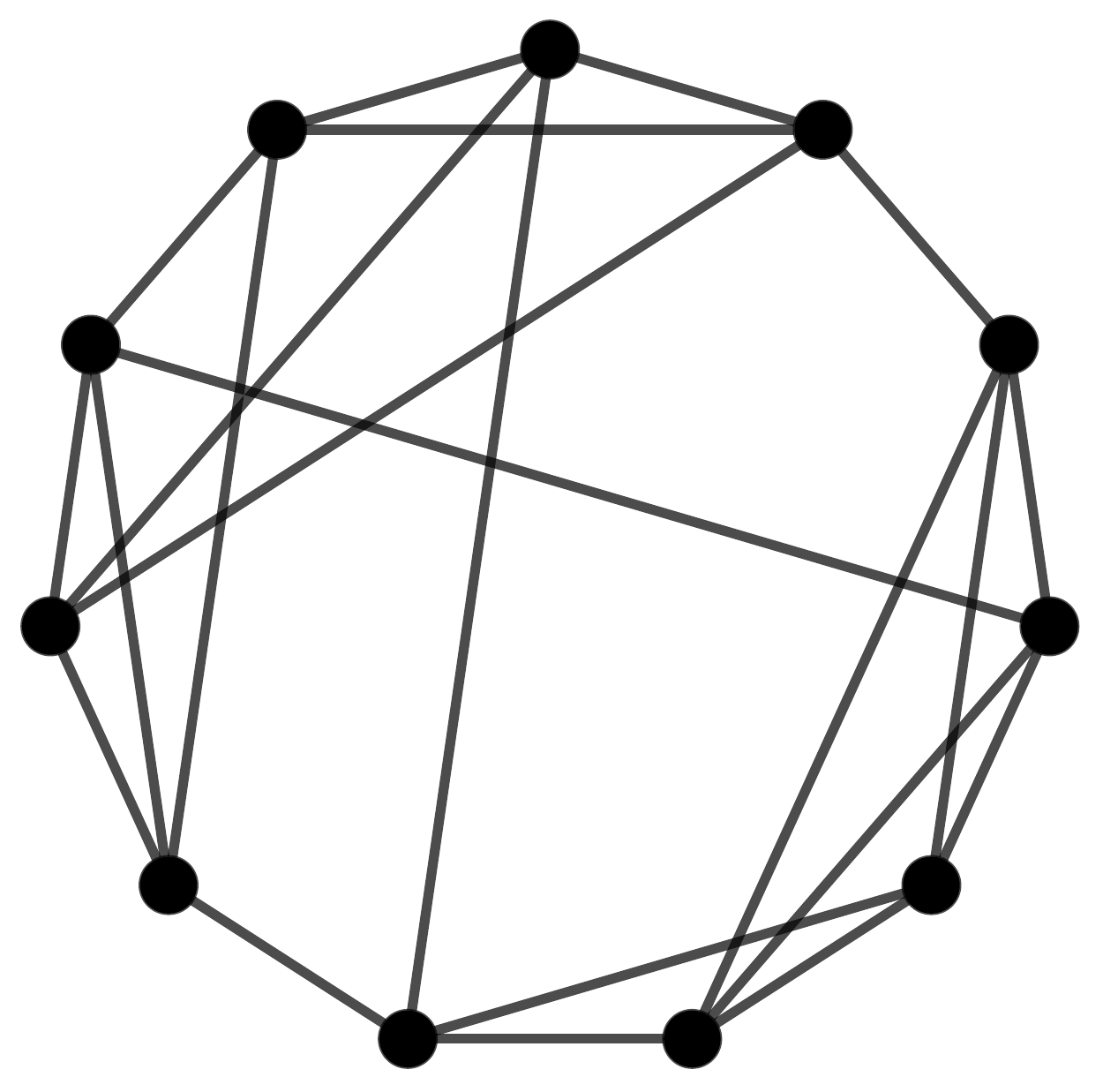}};
\node at (4.5,3.5) {\includegraphics[width=0.34\textwidth]{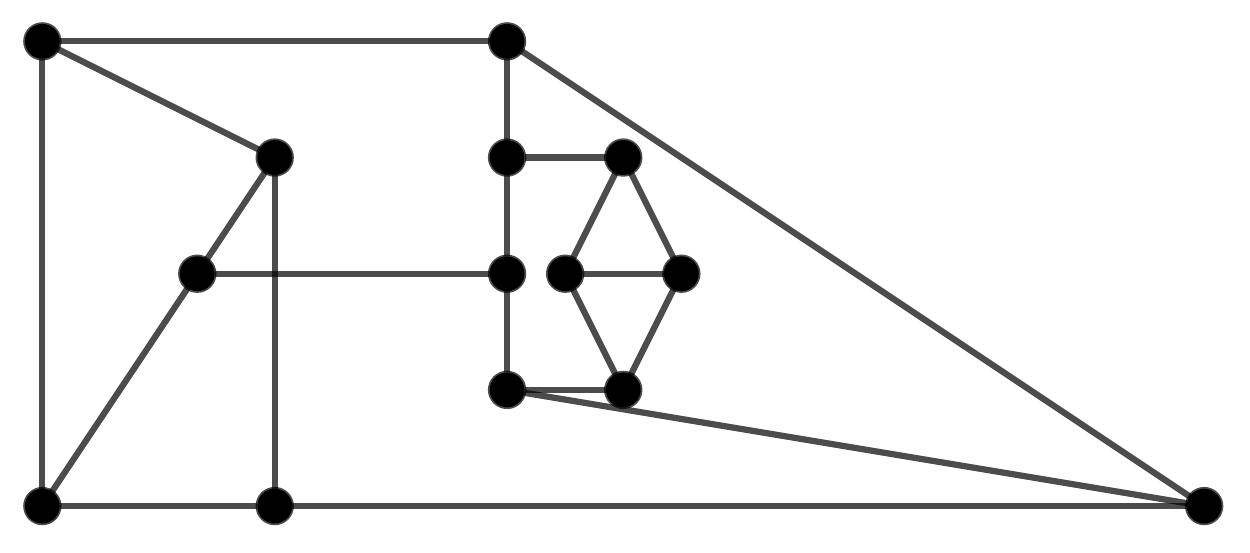}};
\end{tikzpicture}
\vspace{-10pt}
\caption{Two graphs for which $\ker D$ has elements whose entries do not add up to 0. Indeed, $Dw = n \cdot \mathbf{1}$ has no solution.}
\end{figure}
\end{center}

\begin{proof}[Proof of Proposition 3] Let us assume there exist $w_1, w_2 \in \mathbb{R}^n_{\geq 0}$ such that 
$$ Dw_1 = n \cdot \mathbf{1} = Dw_2.$$
Revisiting the proof of Theorem 4, we see that the probability measure
$$ \nu_i = \frac{w_i}{\|w_i\|_{\ell^1}}$$
has the property that for each $u \in V$
$$ \sum_{v \in V}^{} d(u, v) \nu_i(v) = \frac{1}{\|w_i\|_{\ell^1}} \sum_{v \in V}^{} d(u, v)  (w_i)_v = \frac{n}{\|w_i\|_{\ell^1}}.$$
This number then corresponds to the value of $\alpha$ (the value of the game) in the von Neumann Minimax Theorem which is
unique. Therefore, $\|w_1\|_{\ell^1} = \|w_2\|_{\ell^2}$.
\end{proof}

\subsection{Proof of Theorem 3}
\begin{proof} The proof follows quickly from the  standard eigenvalue estimate (see e.g. the textbooks Chung \cite{chung} or Grigor'yan \cite{grig}).
$$ \lambda_1 = \sum_{(u,v) \in E} (f(u)- f(v))^2  \geq \frac{1}{n} \frac{1}{\diam(G)}.$$
We quickly include the very short proof for the convenience of the reader.
Let $f:V \rightarrow \mathbb{R}$ denote an $L^2-$normalized function normalized to
$$ \sum_{v \in V} f(v) = 0$$
so that $f$ minimizes the Dirichlet energy among all $L^2-$normalized functions 
$$ f = \arg \min_{\left\langle f, \mathbf{1} \right\rangle = 0 \atop \|f\|_{ L^2} =1}\sum_{(u,v) \in E} (f(u)- f(v))^2.$$
Since $f$ is $L^2-$normalized, we have 
$$ 1 = \sum_{v \in V} f(v)^2 \leq  \sum_{v \in V} \|f\|_{\ell^{\infty}}^2 = n \cdot \|f\|_{\ell^{\infty}}^2$$
and thus
$\| f\|_{\ell^{\infty}} \geq n^{-1/2}$. Since $f$ has mean value 0, it has to change sign somewhere and therefore
$$ \max_{v \in V} f(v) - \min_{w \in V} f(w) \geq \frac{1}{\sqrt{n}}.$$
At the same time, there is a path $P$ of length at most $\diam(G)$ from the point where the maximum is assumed to
the point where the minimum is assumed. Summing over the path, we get
$$ \frac{1}{\sqrt{n}} \leq \sum_{(u,v) \in P} |f(u) - f(v)| \leq \sqrt{\diam(G)}  \left( \sum_{(u,v) \in P} |f(u) - f(v)|^2 \right)^{1/2}$$
and therefore
$$ \lambda_1 = \sum_{(u,v) \in E} (f(u)- f(v))^2  \geq \frac{1}{n} \frac{1}{\diam(G)}.$$
At this point, we invoke Theorem 1 and conclude that
$$ \lambda_1  \geq \frac{1}{n} \frac{1}{\diam(G)} \geq  \frac{\|w\|_{\ell^1}}{2n^2}.$$

\end{proof}

\subsection{Proof of Theorem 5}

\begin{proof} Suppose that $w \in \mathbb{R}^n_{>0}$ is given. 
There exist two vertices $i,j \in V$ at distance $\diam(G)$ from each other. We abbreviate the average distance between two uniformly at random chosen vertices in the graph as 
$$ \mbox{avdiam}(G) = \frac{1}{n^2} \sum_{u,v \in V} d(u,v) = \frac{1}{n^2} \cdot \left\langle \mathbf{1}, D\mathbf{1}. \right\rangle.$$
The pigeonhole principle implies that there exists a vertex $u \in V$ such that the average distance between $u$ and a uniformly at random chosen vertex $v$ is at least $\mbox{avdiam}(G)$ and therefore
$$ \frac{1}{n} \sum_{v \in V} d(u,v) \geq \mbox{avdiam}(G).$$
Then, inspecting the row of $Dw$ corresponding to the vertex $u$, we see
$$ \| Dw\|_{\ell^{\infty}} \geq   \sum_{v \in V} d(u,v) w_v \geq K \sum_{v \in V} d(u,v) \geq K \cdot n \cdot \mbox{avdiam}(G).$$
This implies
$$ \mbox{avdiam}(G) \leq  \frac{\|Dw\|_{\ell^{\infty}}}{n} \frac{1}{K}.$$
Let us now define the parameter $0 < \delta < 1$ via the equation
 $$\mbox{avdiam}(G) = \delta \diam(G).$$ 
Thus
$$ \diam(G) = \frac{1}{\delta} \mbox{avdiam}(G) \leq  \frac{1}{\delta} \frac{\|Dw\|_{\ell^{\infty}}}{n} \frac{1}{K}$$
Observe that this implies our desired result as soon as $\delta \geq 1/8$. Let now $\delta \leq 1/8$.
If the average diameter is a lot smaller than the diameter, then this implies the existence of a vertex that is fairly close to most other vertices. More precisely, since
$$   \frac{1}{n} \sum_{u \in V}  \frac{1}{n} \sum_{v \in V} d(u,v)  = \mbox{avdiam}(G) = \delta \diam(G)$$
 there has to exist a vertex $u \in V$ such that
$$   \frac{1}{n} \sum_{v \in V} d(u,v) \leq  \delta \diam(G).$$
For this vertex $u$, which will now be fixed for the rest of the proof, there must be many vertices nearby: the set
$$ A = \left\{v \in V: d(u,v) \leq 2 \delta \diam(G) \right\}$$
is necessarily large since
\begin{align*}
 \delta \diam(G) &\geq  \frac{1}{n} \sum_{v \in V} d(u,v)  \geq    \frac{1}{n} \sum_{v \in V \setminus A} d(u,v) \\
 &> \frac{|V \setminus A|}{n} 2 \delta \diam(G) = \frac{n - |A|}{n} 2\delta \diam(G)
\end{align*}
implies that
$$|A| \geq \frac{n}{2}.$$
Let us now pick two vertices $a,b \in V$ at the end-point of a longest path meaning that $d(a,b) = \diam(G)$. Then
$$ \diam(G) = d(a,b) \leq d(a,u) + d(u,b)$$
and thus there exists a vertex $c \in \left\{a,b\right\}$ such that $d(u,c) \geq \diam(G)/2$. Using the triangle inequality
one more time, we see that all the vertices $v \in A$ which are close to $u$ cannot be all that close to $c$ and
$$ \forall~v \in A: \quad d(v,c) \geq \left(\frac{1}{2} - 2\delta\right) \diam(G).$$
Checking now the row of $Dw$ that corresponds to the vertex $c$, we note that all the vertices in $A$ are pretty far away from $c$ and
\begin{align*}
 \| Dw\|_{\ell^{\infty}}  &\geq   \sum_{v \in V} d(c,v) w_v \geq K \sum_{v \in V} d(c,v) \\
 &\geq K \sum_{v \in A} d(c,v) \geq  K \sum_{v \in A} \left(\frac{1}{2} - 2\delta\right) \diam(G) \\
 &=  K |A| \left(\frac{1}{2} - 2\delta\right) \diam(G) \geq K \frac{n}{2} \left(\frac{1}{2} - 2\delta\right) \diam(G).
 \end{align*}
 This now implies
$$ \diam(G) \leq \frac{2}{ \left(\frac{1}{2} - 2\delta\right) n} \frac{\| D w\|_{\ell^{\infty}}}{K}.$$
Recalling that $\delta \leq 1/8$, we arrive at
$$ \diam(G)  \leq  \frac{ \|Dw\|_{\ell^{\infty}}}{n} \frac{ 8}{K}$$
Revisiting the proof of Theorem 3, we deduce from this that
$$ \lambda_1  \geq \frac{1}{n} \frac{1}{\diam(G)} \geq \frac{1}{8 \|Dw \|_{\ell^{\infty}}} K.$$
\end{proof}


\begin{thebibliography}{10}

\bibitem{ambrosio} L. Ambrosio, N.Gigli, G. Savar\'e, 
Metric measure spaces with Riemannian Ricci curvature bounded from below,
Duke Math. J. 163 (2014), no. 7, 1405--1490.

 \bibitem{bakry} D. Bakry and Michel \'Emery, Diffusions hypercontractives, In Seminaire de probabilités XIX 1983/84, pp. 177--206. Springer, Berlin, Heidelberg, 1985.
 
 \bibitem{bauer0} F. Bauer, J. Jost, S. Liu, Ollivier-Ricci curvature and the spectrum of the
normalized graph Laplace operator, Math. Res. Lett. 19.6 (2012), pp. 1185--
1205.

\bibitem{bauer} F. Bauer, F. Chung, Y. Lin, Y. Liu, Curvature aspects of graphs.  Proc. Amer. Math. Soc 145 (2017), p. 2033--2042

\bibitem{bauer2} F. Bauer, B. Hua, J. Jost, S. Liu and G. Wang, The geometric meaning of
curvature: local and nonlocal aspects of Ricci curvature, in: `Modern approaches
to discrete curvature' Lecture Notes in Math. 2184, p. 1--62, Springer, Cham,
2017.

\bibitem{bourne} D. Bourne, D. Cushing, S. Liu, F. M\"unch, and N. Peyerimhoff. Ollivier--Ricci Idleness Functions of Graphs. SIAM Journal on Discrete Mathematics 32, no. 2 (2018): 1408--1424.

\bibitem{cheng} S. Y. Cheng, Eigenvalue comparison theorems and its geometric applications, Mathematische Zeitschrift, 143 (1975): p. 289--297

\bibitem{chung} F. Chung, Spectral Graph Theory, CBMS Regional Conference Series in Mathematics 92, American Mathematical Society.

\bibitem{cleary} J. Cleary and S. A. Morris, Numerical geometry-numbers for shapes, Amer. Math. Monthly 93 (1986) 260--275.

\bibitem{cush1} D. Cushing, R. Kangaslampi, Y. Lin, S. Liu, L. Lu, and S.-T. Yau, Erratum for Ricci-flat graphs with girth at least five, to appear in Communications in Analysis and Geometry, arXiv 1802:02979;

\bibitem{cush2} D. Cushing, R. Kangaslampi, Y. Lin, S. Liu, L. Lu, and S.-T. Yau, Ricci-flat cubic graphs with girth five. arXiv preprint arXiv:1802.02982 (2018).

\bibitem{rigidity} D. Cushing, S. Kamtue, J. Koolen, S.Liu, F. M\"unch, N.Peyerimhoff,  Rigidity of the Bonnet-Myers inequality for graphs with respect to Ollivier Ricci curvature, Advances in Mathematics 369 (2020), 107--188

\bibitem{gross} O. Gross, The rendezvous value of a metric space, in: Advances in Game Theory, Ann. of Math
 Studies no. 52, Princeton (1964) 49-53.

\bibitem{forman} Robin Forman, Bochner’s method for cell complexes and combinatorial Ricci curvature. Discrete
and Computational Geometry, 29(3): 323--374, 2003.

\bibitem{glicksberg} I. Glicksberg, A Further Generalization of the Kakutani Fixed Point Theorem, with Application to Nash Equilibrium Points, Proceedings of the American Mathematical Society 3 (1952), p. 170--174

\bibitem{grig} A. Grigor’yan, Introduction to analysis on graphs, University Lecture Series 71, American Mathematical Society, 2018

\bibitem{higuchi}  Y. Higuchi, Combinatorial curvature for planar graphs, J. Graph Theory 38
(2001), 220--229.

\bibitem{horn} P. Horn, Y. Lin, S. Liu, S.-T. Yau, Volume doubling, Poincar\'e inequality
and Gaussian heat kernel estimate for non-negatively curved graphs, J. Reine
Angew. Math. 757 (2019), 89–130.

\bibitem{jost} J. Jost, S. Liu, Ollivier’s Ricci Curvature, Local Clustering
and Curvature-Dimension Inequalities on Graphs, Discrete Comput. Geom. 51 (2014), p. 300--322

\bibitem{lich}  A. Lichnerowicz,  G\'eom\'etrie des groupes de transformations, Dunod, Paris, 1958.

\bibitem{lly} Y. Lin, L. Lu, S.-T. Yau, Ricci curvature of graphs, Tohoku Mathematical Journal, Second Series 63, no. 4 (2011): 605--627.

\bibitem{lly2} Y. Lin, L. Lu, S.-T. Yau, Ricci-flat graphs with girth at least five, Comm. Anal. Geom. 22 (2014), no. 4, 671--687.

\bibitem{linyau} Y. Lin and S.-T. Yau, 
Ricci curvature and eigenvalue estimate on locally finite graphs. 
Math. Res. Lett. 17 (2010), no. 2, 343--356.

\bibitem{lott} J. Lott and C. Villani, Ricci curvature for metric-measure spaces via optimal transport, Ann. of Math. (2)
169 (2009), p. 903--991.

\bibitem{maas} J. Maas, Gradient flows of the entropy for finite Markov chains,
J. Funct. Anal. 261 (2011), no. 8, p. 2250--2292.

\bibitem{myers} S. B. Myers, Riemannian manifolds with positive mean curvature, Duke Mathematical Journal, 8 (1941): p.401--404

\bibitem{olli0} Y. Ollivier, Ricci curvature of metric spaces, C. R. Math. Sci. 345.11 (2007), pp. 643--646.

\bibitem{ollivier} Y. Ollivier,
Ricci curvature of Markov chains on metric space,
J. Funct. Anal., 256 (2009), pp. 810-864

\bibitem{olli2} Y. Ollivier, A survey of Ricci curvature for metric spaces and Markov chains. Probabilistic approach to geometry, 343--381, Adv. Stud. Pure Math., 57, Math. Soc. Japan, Tokyo, 2010. 

\bibitem{villani} Y. Ollivier and C. Villani, A curved Brunn-Minkowski inequality on the discrete
hypercube, or: What is the Ricci curvature of the discrete hypercube?, SIAM J.
Discrete Math. 26(3) (2012), 983--996.

\bibitem{stein} S. Steinerberger, The first eigenvector of a distance matrix is nearly constant, arXiv:2205.15920

\bibitem{stone} D. A. Stone, A combinatorial analogue of a theorem of Myers, Illinois J. Math.
20(1) (1976), p. 12--21 and Correction to my paper: A combinatorial analogue of
a theorem of Myers, Illinois J. Math. 20(3) (1976), 551--554.

\bibitem{sturm}  K.-T. Sturm, On the geometry of metric measure spaces, (I), (II), Acta Math. 196 (2006), 65--131, 133--177.

\bibitem{thom} C. Thomassen, The rendezvous number of a symmetric matrix and a compact connected metric space. Amer. Math. Monthly 107 (2000), no. 2, 163--166. 

\bibitem{john} J. von Neumann, Zur Theorie der Gesellschaftsspiele, Math. Ann. 100 (1928): p. 295--320. 

\bibitem{woess} W. Woess, A note on tilings and strong isoperimetric inequality, Math. Proc.
Cambridge Philos. Soc. 124(3) (1998), p. 385--393.

\end{thebibliography}
\end{document}